\documentclass[11pt]{article}

\usepackage{latexsym,amsmath,amsfonts,amssymb,amstext,theorem,euscript,
array,enumerate,mathrsfs}
\usepackage{color}
\usepackage{comment}

\newtheorem{Theorem}{Theorem}[section]

\newtheorem{Proposition}{Proposition}[section]

\newtheorem{Lemma}{Lemma}[section]

\newtheorem{Remark}{Remark}[section]

\numberwithin{equation}{section}

\def\esssup_#1{\underset{#1}{\mathrm{ess\,sup\, }}}
\def\essinf_#1{\underset{#1}{\mathrm{ess\,inf\, }}}

\def \trans{^{\scriptscriptstyle{\intercal}}}

\def\sqr#1#2{{\vcenter{\vbox{\hrule height .#2pt \hbox{\vrule
 width .#2pt height#1pt \kern#1pt \vrule
width .#2pt} \hrule height .#2pt}}}}

\def\ds{\begin{displaystyle}}
\def\eds{\end{displaystyle}}
\def\dis{\displaystyle }
\def\<{\langle }
\def\>{\rangle }

\def \N{\mathbb{N}}
\def \R{\mathbb{R}}

\def \E{\mathbb{E}}
\def \F{\mathbb{F}}

\def \P{\mathbb{P}}

\def \Ac{{\cal A}}
\def \Bc{{\cal B}}
\def \Cc{{\cal C}}

\def \Ec{{\cal E}}
\def \Fc{{\cal F}}

\def \Lc{{\cal L}}
\def \Pc{{\cal P}}
\def \Mc{{\cal M}}

\def \Rc{{\cal R}}

\def \Vc{{\cal V}}

\def \Vc{{\cal V}}

\def\calf{{\cal F}}

\def\bfC{{\bf C}}
\def\bfM{{\bf M}}

\def \eps{\varepsilon}

\def \ep{\hbox{ }\hfill$\Box$}

\allowdisplaybreaks

\begin{document}

\title{Backward SDEs and infinite horizon stochastic optimal control}

\author{Fulvia CONFORTOLA\thanks{Politecnico di Milano, Dipartimento di Matematica, via Bonardi 9, 20133 Milano, Italy; e-mail: \texttt{fulvia.confortola@polimi.it}}
\and
Andrea COSSO\thanks{Politecnico di Milano, Dipartimento di Matematica, via Bonardi 9, 20133 Milano, Italy; e-mail: \texttt{andrea.cosso@polimi.it}}
\and
Marco FUHRMAN\thanks{Universit\`a di Milano, Dipartimento di Matematica, via Saldini 50, 20133 Milano,
Italy; e-mail: \texttt{marco.fuhrman@unimi.it} This author was supported by the Italian
MIUR-PRIN 2015 ``Deterministic and stochastic evolution equations''.}
}

\maketitle

\begin{abstract}
We study an optimal control problem on infinite horizon for a controlled
stochastic differential equation driven by Brownian motion, with a discounted
reward functional. The equation may have memory or delay effects
in the coefficients, both with respect to state and control,
and the noise can be degenerate.
We prove that the value, i.e. the supremum
of the reward functional over all admissible controls, can be represented
by the solution of an associated backward stochastic
differential equation (BSDE) driven by the Brownian motion and
an auxiliary independent Poisson process and having a sign constraint
on jumps.

In the Markovian case when the coefficients
depend only on the present values of the state and the control,
we prove that the BSDE can be used
to construct the solution, in the sense of viscosity theory,
to the corresponding Hamilton-Jacobi-Bellman partial
differential equation of elliptic type on the whole space, so that
it provides us with a Feynman-Kac representation in this
fully nonlinear context.

The method of proof consists in showing that the value
of the original problem
is the same as the value
of an auxiliary optimal control problem (called randomized), where
the control process is replaced by a fixed pure jump process
and maximization is taken over a class of absolutely continuous changes
of measures which affect the stochastic intensity of the
jump process but leave the law of the driving Brownian motion unchanged.
\end{abstract}

\vspace{5mm}


\vspace{5mm}


\maketitle

\date{}

\newpage

\section{Introduction}

Let us consider a classical optimal control problem with infinite horizon for a
stochastic  equation
in $\R^n$
of the form
\begin{equation}\label{sdeintro}
X_t^\alpha \ = \ x + \int_0^t b(X^\alpha_s,\alpha_s)\,ds + \int_0^t \sigma(X_s^\alpha,\alpha_s)\,dW_s,
 \qquad t\geq0,
\end{equation}
starting at a point $x\in\R^n$,
with discounted reward functional
\[
J(x,\alpha) \ = \ \E\bigg[\int_0^\infty e^{-\beta t} f(X_t^\alpha,\alpha_t)\,dt\bigg].
\]
Here a
$W$ is a $d$-dimensional Brownian motion defined in some probability space,
$b$ and $\sigma$ are given coefficients with values in $\R^n$
and $\R^{n\times d}$ respectively,
$f$ is a real function representing the running cost rate, $\beta>0$
is a discount factor. The control
process
$\alpha$ is a stochastic process taking values in a
metric space $A$ and progressive with respect to the
completed filtration $\F^W$ generated by the Brownian motion;
the class of such controls is denoted by $\Ac$.
The aim is to maximize $J(x,\alpha)$ over $\Ac$
and to characterize the value function
$$
v(x)=\sup_{\alpha \in\Ac}
J(x,\alpha).
$$

It is well known that, under natural assumptions, the value function
is well defined and it is a viscosity solution to the Hamilton-Jacobi-Bellman
(HJB)
equation, which is the following  elliptic partial differential equation
on the whole space:
\begin{equation}\label{HJBintro}
\beta\,v(x) - \sup_{a\in A}\big[\Lc^a v(x) + f(x,a)\big] \ = \ 0, \qquad x\in\R^n,
\end{equation}
where $\Lc^a$ is the
Kolmogorov operator depending on the control parameter $a$:
\[
\Lc^a v(x) \ = \ \langle b(x,a),D_x v(x)\rangle + \frac{1}{2}\text{tr}\big[\sigma(x,a)
\sigma\trans(x,a)D_x^2 v(x)\big].
\]
 Thus, when a uniqueness result holds, the  HJB equation completely characterizes
the value function.

It is the purpose of this paper to provide a different
representation of the value function, based on backward stochastic
differential equations (BSDEs).
BSDEs are used since long to represent value functions of stochastic optimal
control problems and more generally solutions to partial differential
equations of parabolic and elliptic type.
Besides their intrinsic interest, some motivations are the fact that they
usually allow to extend the obtained results beyond the Markovian case and
they often admit  efficient numerical approximations, which makes them a
competitive tool in comparison with other more common numerical
methods for partial differential equations.
However the classical results relating BSDEs with partial differential equations,
see for instance \cite{Peng91}  and \cite{PardouxPeng92},
only cover  cases when the equation is semilinear,
i.e. it takes the form (in the elliptic case)
\begin{equation}\label{PDEintro}
\beta\,v(x) -  \Lc v(x) + \psi(x,v(x),D_xv(x)\sigma(x)) \ = \ 0, \qquad x\in\R^n,
\end{equation}
where $\Lc$ is the linear Kolmogorov operator associated to (uncontrolled)
coefficients $b(x)$, $\sigma(x)$, and the nonlinear term $\psi$ depends on the
gradient $D_xv(x)$ only through the product $D_xv(x)\sigma(x)$.
The general HJB equation \eqref{HJBintro}, being fully nonlinear, can not be cast in the
form \eqref{PDEintro} except in special cases. The corresponding
optimal control problems also have a special form,
and in particular the occurrence of a diffusion coefficient
$\sigma(x,a)$ depending on a control parameter $a\in A$ can not be allowed.

To overcome this difficulty,  new methods have recently been developed.
We mention the theory of second order BSDEs  \cite{STZ}
and the theory of $G$-expectations \cite{Peng07}, that both allow a probabilistic
representation of solutions to classes of  fully nonlinear equations.

In this paper we follow a different approach, based on a method
that we call {\it randomization of control}. It has been introduced
in \cite{bou09} and then successfully
applied to several stochastic optimization problems,
including impulse control, optimal switching, optimal stopping:
see \cite{EKa}, \cite{EKb},
\cite{KMPZ10},
\cite{FuhrmanPhamZeni16}, \cite{BandiniFuhrman15},
\cite{Bandini15} and especially \cite{KP12} for
a systematic application to a generalized class of HJB equations.
In the proofs of our result we will especially follow
\cite{FP15} and \cite{BCFP16b}
that deal with the optimal control problem with finite horizon
and, in the Markovian case, with the corresponding parabolic HJB equation.

The general idea of the randomization method is as follows. By enlarging the original probability space
if necessary, we consider an independent Poisson random measure
$\mu(dt,da)$ on $(0,\infty)\times A$, with finite intensity measure
$\lambda(da)$, and the corresponding
$A$-valued process with piecewise constant trajectories, denoted by $I$.
We formally replace the control process $\alpha$ by $I$, so we
solve the equation
\begin{equation}\label{sdeintrorandomizzata}
X_t \ = \ x + \int_0^t b(X_s,I_s)\,ds + \int_0^t \sigma(X_s,I_s)\,dW_s,
 \qquad t\geq0.
\end{equation}
Then we consider an auxiliary optimization problem, called randomized problem,
which consists in optimizing among equivalent changes of probability measures which only
affect the intensity measure of $I$ but not the law of $W$.
In the randomized problem, an admissible
control is a bounded positive map $\nu$ defined on
$\Omega\times (0,\infty) \times A$, which is predictable with respect to
the filtration $\F^{W,\mu}$ generated by $W$ and $\mu$.
Given $\nu$, by means of an absolutely continuous change
of probability measure of Girsanov type we construct a probability $\P^\nu$
such that the compensator of
$\mu$ is given by $\nu_t(a)\lambda(da)dt$ and $W$ remains a Brownian motion under
$\P^\nu$. Then we introduce another
reward functional and the corresponding value function
\[
J^\Rc(x,\nu) \ = \ \E^\nu\bigg[\int_0^\infty e^{-\beta t} f(X_t,I_t)\,dt\bigg],
\qquad
v^\Rc(x)=\sup_{\nu}
J^\Rc(x,\nu),
\]
where $\E^\nu$ denotes the expectation under $\P^\nu$.
Some technical issues arise in connection with the use of Girsanov
transformation on the whole time halfline, so the precise formulation
is slightly more involved: see section
\ref{S:Randomized} below and especially Remark \ref{girsanovinfinite}
for more details.
Building on our previous results in \cite{FP15} and \cite{BCFP16b}
we prove that the two value functions $v$ and $v^\Rc$
coincide: see  Theorem \ref{T:Identification}.

Next we prove that the function $v^\Rc$ can be represented
by means of the following infinite horizon backward
stochastic differential equation (BSDE) with constraint: for all $0 \leq t \leq T<\infty$,
\begin{equation} \label{BSDE_Markov_intro}
\left\{
\begin{array}{l}
Y_t   =  \dis  Y_T - \beta\int_t^T Y_s\,ds + \int_t^T f( X_s, I_s)\,ds + K_T - K_t
 - \int_t^T Z_s\,d W_s - \int_t^T \int_A U_s(a)\,\mu(ds,da),
 \\
U_t(a)   \leq  0.
\end{array}\right.
\end{equation}
In this equation the unknown process is a quadruple $(Y_t,Z_t,U_t(a),K_t)$
where $Y$ is c\`adl\`ag adapted (with respect to $\F^{W,\mu}$),
$Z$ is $\R^d$-valued progressive, $K$ is increasing predictable,
$U$ is a predictable random field, all
satisfying appropriate integrability conditions. The constraint on $U$ can be seen as
a constraint of
nonpositivity on the totally inaccessible jumps of $Y$.
We prove that the BSDE \eqref{BSDE_Markov_intro} admits a unique
minimal solution in a suitable sense and that it represents the value function
in the sense that $Y_0=v^\Rc(x)$ and so also $Y_0=v(x)$, so that
 we obtain the desired representation of the value function
for our original control problem (see Theorem \ref{T:Main}).
In addition, we prove that the solution to the BSDE satisfies
a sort of recursive formula (formula
\eqref{DPP}) which is a version of the dynamic programming principle
in the setting of the randomized control problem. We exploit this
functional equality to prove that the value function $v$ is a viscosity
solution to the HJB equation \eqref{HJBintro}, see
Theorem \ref{T:Markov}. Therefore the equality $Y_0=v(x)$ can also
be seen as a fully nonlinear Feynman-Kac representation
for the solution to \eqref{HJBintro}.
This approach allows to circumvent the difficulties related to a rigorous
proof of the classical dynamic programming principle which might
be lengthy and, in some versions, may require the use of nontrivial
measurability arguments. On the contrary,
we do not deal with uniqueness results
for the HJB equation, which are classical and are known to hold under suitable assumptions,
see Remark \ref{HJBunique}. We also stress that none of our results
requires nondegeneracy assumption on the noise, so no
requirements are imposed on the diffusion coefficient $\sigma$
except Lipschitz conditions and some  boundedness and continuity assumptions.

As mentioned before, the use  of BSDEs often allows for
efficient numerical treatment. This is also the case for the constrained
BSDEs of the form \eqref{BSDE_Markov_intro}, at least in the finite horizon
case: see \cite{KLP14}, \cite{KLP15}.

Another advantage of our technique is that we are able to generalize all the previous
results (except the ones on the HJB equation) to the non Markovian case when the
coefficients of the controlled equation exhibit memory effects, i.e.
for an equation of the form
\begin{equation}\label{sde_intro}
X_t^\alpha \ = \ x + \int_0^t b_s(X^\alpha,\alpha)\,ds + \int_0^t \sigma_s(X^\alpha,\alpha)\,dW_s,
\qquad t\ge 0,
\end{equation}
where, at any time $t$, the value
of the coefficients $b_t(X^\alpha,\alpha)$ and $\sigma_t(X^\alpha,\alpha)$
may depend on the entire past trajectory of the state $(X^\alpha_s)_{s\in [0,t]}$
as well as the control process $(\alpha_s)_{s\in [0,t]}$.
In fact, the previous results are formulated and proved directly in this generality,
while the Markovian case is only addressed to deal with the HBJ equation.

We finally mention that in the paper
\cite{CossoFuhrmanPham16}, co-authored by some of us,
 BSDEs of the form \eqref{BSDE_Markov_intro} have been introduced,
as an intermediate technical step in the proofs, when dealing
with HJB equations of ergodic type, namely of the form
$$
\lambda - \sup_{a\in A}\big[\Lc^a v(x) + f(x,a)\big] \ = \ 0, \qquad x\in\R^n,
$$
where both the function $v$ and the constant $\lambda$ are unknown.
However, the results in \cite{CossoFuhrmanPham16} impose strong
restrictions on the coefficients $b,\sigma,f$ and the space of control
actions $A$. In particular, Lipschitz conditions
were imposed on $f$ and special dissipativity assumptions
were imposed on $b$ and $\sigma$ in order to guarantee appropriate
ergodicity properties. In the present paper these assumptions are  dropped. In addition,
the results in \cite{CossoFuhrmanPham16} depend in an essential
way on the Markovianity of the stochastic system and can not
be applied to the controlled equation \eqref{sde_intro}.  As a result, we are led
to a careful study of the growth rate of the solution $X$ to
the non Markovian equation \eqref{sde_intro}
(compare Lemma \ref{L:Estimate} below)
and we have to relate it to polynomial growth conditions imposed
on $f$ as well as an appropriate value for the discount factor $\beta$.

The plan of the paper is as follows.
In section \ref{S:Formulation} we formulate our assumptions
for the general non Markovian framework and introduce
the optimal control problem on infinite horizon, with special
attention to the behaviour of the controlled system for large times, whereas
in section \ref{S:Randomized} we formulate the auxiliary randomized problem.
In section \ref{S:BSDE} we prove the equality of the values
of these two problems, we introduce and study the well-posedness
of the  constrained BSDE \eqref{BSDE_Markov_intro} and we
prove that it gives the desired representation of the values.
Finally, in section \ref{S:Markov}, we restrict to the Markovian case
and prove that the solution to the BSDE provides us with
a viscosity solution to the HJB equation \eqref{HJBintro}.

\section{Formulation of the infinite horizon optimal control problem}
\label{S:Formulation}

\setcounter{equation}{0} \setcounter{Assumption}{0}
\setcounter{Theorem}{0} \setcounter{Proposition}{0}
\setcounter{Corollary}{0} \setcounter{Lemma}{0}
\setcounter{Definition}{0} \setcounter{Remark}{0}

Let $(\Omega,\Fc,\P)$ be a complete probability space, on which a $d$-dimensional Brownian motion $W=(W_t)_{t\geq0}$ is defined. Let $\F^W=(\Fc_t^W)_{t\geq0}$ denote the $\P$-completion of the filtration generated by $W$. Let $A$ be a nonempty Borel space (namely, $A$ is a topological space homeomorphic to a Borel subset of a Polish space) and denote by $\Ac$ the set of $\F^W$-progressive processes $\alpha\colon\Omega\times[0,\infty)\rightarrow A$. $A$ is the space of control actions and $\Ac$ is the family of admissible control processes. Finally, we denote $\Bc(A)$ the Borel $\sigma$-algebra of $A$.

Fix a deterministic point $x_0\in\R^n$. For every $\alpha\in\Ac$, consider the controlled equation:
\begin{equation}\label{sde}
X_t^\alpha \ = \ x_0 + \int_0^t b_s(X^\alpha,\alpha)\,ds + \int_0^t \sigma_s(X^\alpha,\alpha)\,dW_s,
\end{equation}
for all $t\geq0$. The infinite horizon stochastic optimal control problem consists in maximixing over $\alpha\in\Ac$ the gain functional
\[
J(\alpha) \ = \ \E\bigg[\int_0^\infty e^{-\beta t} f_t(X^\alpha,\alpha)\,dt\bigg].
\]
The constant $\beta>0$ will be specified later. The coefficients $b$, $\sigma$, $f$ are defined on $[0,\infty)\times \bfC_n\times  {\bf M}_A$ with
values in $\R^n$, $\R^{n\times d}$, $\R$, respectively, where:
\begin{itemize}
\item $\mathbf C_n$ is the set of continuous trajectories $x\colon[0,\infty)\rightarrow\R^n$. We introduce the canonical filtration $(\mathcal C_t^n)_{t\geq0}$ and denote $Prog(\mathbf C_n )$ the $(\Cc_t^n)$-progressive $\sigma$-algebra on $[0,\infty)\times\mathbf C_n $;
\item $\mathbf M_A$ is the  set  of Borel measurable trajectories $a\colon[0,\infty)\rightarrow A$. We introduce the canonical filtration
$(\Mc_t^A)_{t\in[ 0,T]}$ and denote $Prog(\mathbf C_n \times \mathbf M_A)$ the $(\Cc_t^n\otimes \Mc_t^A)$-progressive $\sigma$-algebra on $[0,\infty)\times\mathbf C_n\times\mathbf M_A$.
\end{itemize}

In the present paper, we consider the two following alternative sets of assumptions on $b$, $\sigma$, $f$. Notice that {\bf (A)} differs from {\bf (A)'} only for points (iii) and (iv).

\vspace{3mm}

\noindent {\bf (A)}
\begin{itemize}
\item [(i)]
The functions $b$, $\sigma$, $f$ are $Prog(\mathbf C_n \times \mathbf M_A)$-measurable.
\item[(ii)] For every $T>0$, if $x_m,x\in\bfC_n$, $a_m,a\in
\bfM_A$, $\sup_{t\in[0,T]}|x_m(t)-x(t)|\to 0$, $a_m(t)\to a(t)$ for $dt$-a.e.
  $t$ $\in$ $[0,T]$ as $m\to\infty$, then we have
  $$
b_t(x_m,a_m) \to b_t(x ,a),\;\;  \sigma_t(x_m,a_m)\to \sigma_t(x,a),\;\;  f_t(x_m,a_m)\to f_t(x,a),
\quad
\text{ for } dt\text{-a.e. }
 t \in [0,T].
 $$
\item [(iii)]
For every $T>0$, there exists a constant  $L_T\geq0$ such that
\begin{align*}
|b_t(x,a) - b_t(x',a)| + |\sigma_t(x,a)-\sigma_t(x',a)|
\ &\leq \ L_T \sup_{s\in[0,t]}|x(s)-x'(s)|,
\\
|b_t(0,a)| + |\sigma_t(0,a)| \ &\leq \ L_T,
\end{align*}
for all $t\in[0,T]$, $x,x'\in\mathbf C^n$, $a\in\mathbf M_A$.
\item[(iv)] The function $f$ is bounded. We denote $\|f\|_\infty:=\sup_{t,x,a}|f_t(x,a)|<\infty$.
\item [(v)] $\beta$ can be any strictly positive real number.
\end{itemize}

\vspace{2mm}

\noindent {\bf (A)'}
\begin{itemize}
\item [(i)]
The functions $b$, $\sigma$, $f$ are $Prog(\mathbf C_n \times \mathbf M_A)$-measurable.
\item[(ii)] For every $T>0$, if $x_m,x\in\bfC_n$, $a_m,a\in
\bfM_A$, $\sup_{t\in[0,T]}|x_m(t)-x(t)|\to 0$, $a_m(t)\to a(t)$ for $dt$-a.e.
  $t$ $\in$ $[0,T]$ as $m\to\infty$, then we have
  $$
b_t(x_m,a_m) \to b_t(x ,a),\;\;  \sigma_t(x_m,a_m)\to \sigma_t(x,a),\;\;  f_t(x_m,a_m)\to f_t(x,a),
\quad
\text{ for } dt\text{-a.e. }
 t \in [0,T].
 $$
\item [(iii)]
There exists a constant  $L\geq0$ such that
\begin{align*}
|b_t(x,a) - b_t(x',a)| + |\sigma_t(x,a)-\sigma_t(x',a)|
\ &\leq \ L \sup_{s\in[0,t]}|x(s)-x'(s)|,
\\
|b_t(0,a)| + |\sigma_t(0,a)| \ &\leq \ L,
\end{align*}
for all $t\geq0$, $x,x'\in\mathbf C^n$, $a\in\mathbf M_A$.
\item[(iv)] There exist constants $M>0$ and $r>0$ such that
\[
|f_t(x,a)| \ \leq \ M(1 + \sup_{s\in[0,t]}|x(s)|^r),
\]
for all $t\geq0$, $x\in\mathbf C^n$, $a\in\mathbf M_A$.
\item[(v)] $\beta>\bar\beta$, where $\bar\beta$ is a strictly positive real number such that
\begin{equation}\label{beta}
\E\Big[\sup_{s\in[0,t]}|X_s^\alpha|^r\Big] \ \leq \ \bar C e^{\bar\beta t}(1 + |x_0|^r),
\end{equation}
for some constant $\bar C\geq0$, with $\bar\beta$ and $\bar C$ independent of $t\geq0$, $\alpha\in\Ac$, $x_0\in\R^n$. See Lemma \ref{L:Estimate}.
\end{itemize}

Notice that under either {\bf (A)}-(i)-(ii)-(iii) or {\bf (A)'}-(i)-(ii)-(iii), there exists a unique $\F^W$-progressive continuous process $X^\alpha=(X_t^\alpha)_{t\geq0}$ solution to equation \eqref{sde}. Moreover, under {\bf (A)}-(i)-(ii)-(iii), for every $T>0$ and $p>0$, we have that there exists a constant $C_{T,p}\geq0$ such that
\[
\E\Big[\sup_{t\in[0,T]}|X_t^\alpha|^p\Big] \ \leq \ C_{T,p}\big(1 + |x_0|^p\big).
\]
On the other hand, under {\bf (A)'}-(i)-(ii)-(iii), we have that there exists $\bar\beta\geq0$ such that estimate \eqref{beta} holds. This latter estimate follows from the next Lemma \ref{L:Estimate}, where we prove a more general result needed later.
\begin{Lemma}\label{L:Estimate}
Suppose that Assumption {\bf (A')}\textup{-(i)-(ii)-(iii)} holds. Let $(\hat\Omega,\hat\Fc,\hat\P)$ be a complete probability space and let $\hat\F=(\hat\Fc_t)_{t\geq0}$ be a filtration satisfying the usual conditions. Let also $\hat W=(\hat W_t)_{t\geq0}$ be a $d$-dimensional Brownian motion on the filtered space $(\hat\Omega,\hat\Fc,\hat\F,\hat\P)$. Let $\gamma\colon\hat\Omega\times[0,\infty)\rightarrow A$ be an $\hat\F$-progressive process. Finally, let $\hat X=(\hat X_t)_{t\geq0}$ be the unique continuous $\hat\F$-adapted process solution to the following equation
\begin{equation}\label{sde_hat}
\hat X_t \ = \ x_0 + \int_0^t b_s(\hat X,\gamma)\,ds + \int_0^t \sigma_s(\hat X,\gamma)\,d\hat W_s,
\end{equation}
for all $t\geq0$. Then, for every $p>0$, there exist two constants $\bar C_{p,L}\geq0$ and $\bar\beta_{p,L}\geq0$ such that
\begin{equation}\label{beta_hat}
\hat\E\Big[\sup_{s\in[0,T]}|\hat X_s|^p\Big|\hat\Fc_t\Big] \ \leq \ \bar C_{p,L}\,e^{\bar\beta_{p,L}\,(T-t)}\,(1 + \sup_{s\in[0,t]}|\hat X_s|^p), \qquad \hat\P\text{-a.s.}
\end{equation}
for all $t\geq0$, $T\geq t$, with $\bar C_{p,L}$ and $\bar\beta_{p,L}$ depending only on $p$ and the constant $L$ appearing in Assumption {\bf (A')}\textup{-(iii)}. When $p=r$, with $r$ as in Assumption {\bf (A')}\textup{-(iv)}, we denote $\bar C_{r,L}$ and $\bar\beta_{r,L}$ simply by $\bar C$ and $\bar\beta$.
\end{Lemma}
\textbf{Proof}
See Appendix. We remark that, to our knowledge, a proof of estimate \eqref{beta_hat} in the path-dependent case does not exist in the literature. On the other hand, for the non-path-dependent case we refer for instance to Theorem II.5.9 in \cite{80Krylov}. Notice however that those proofs use in an essential way the fact that $b$ and $\sigma$ depends only on the present value of the process $\hat X$, so that they can not be extended to the path-dependent case.
\ep

\vspace{3mm}

We define the value of the control problem as
\[
V \ = \ \sup_{\alpha\in\Ac} J(\alpha).
\]
Now, for every $T>0$, consider the finite horizon optimal control problem with value
\[
V_T \ = \ \sup_{\alpha\in\Ac} J_T(\alpha),
\]
where
\[
J_T(\alpha) \ = \ \E\bigg[\int_0^T e^{-\beta t}\,f_t(X^\alpha,\alpha)\,dt\bigg].
\]
\begin{Lemma}\label{L:ApproxV}
Under either {\bf (A)} or {\bf (A)'}, we have
\[
V \ = \ \lim_{T\rightarrow\infty} V_T.
\]
\end{Lemma}
\textbf{Proof.}
\emph{Assumption {\bf (A)} holds.} We have
\[
|V - V_T| \ \leq \ \sup_{\alpha\in\Ac} \E\bigg[\int_T^\infty e^{-\beta t} \big|f_t(X^\alpha,\alpha)\big|\,dt\bigg] \ \leq \ \|f\|_\infty \frac{e^{-\beta T}}{\beta} \ \overset{T\rightarrow\infty}{\longrightarrow} \ 0.
\]

\vspace{1mm}

\noindent\emph{Assumption {\bf (A)'} holds.} We have
\begin{align*}
|V - V_T| \ \leq \ \sup_{\alpha\in\Ac} \E\bigg[\int_T^\infty e^{-\beta t} \big|f_t(X^\alpha,\alpha)\big|\,dt\bigg] \ \leq \ M \frac{e^{-\beta T}}{\beta} + M\E\bigg[\int_T^\infty e^{-\beta t} \sup_{s\in[0,t]}|X_s^\alpha|^r\,dt\bigg]& \\
\leq \ M \frac{e^{-\beta T}}{\beta} + M\int_T^\infty e^{-(\beta-\bar\beta)t} \E\Big[e^{-\bar\beta t}\sup_{s\in[0,t]}|X_s^\alpha|^r\Big]\,dt&.
\end{align*}
Using \eqref{beta}, we obtain
\[
|V - V_T| \ \leq \ M \frac{e^{-\beta T}}{\beta} + M \bar C \big(1 + |x_0|^r\big) \frac{e^{-(\beta-\bar\beta)T}}{\beta - \bar\beta} \ \overset{T\rightarrow\infty}{\longrightarrow} \ 0.
\]
\ep

\section{Randomized optimal control problem}
\label{S:Randomized}

\setcounter{equation}{0} \setcounter{Assumption}{0}
\setcounter{Theorem}{0} \setcounter{Proposition}{0}
\setcounter{Corollary}{0} \setcounter{Lemma}{0}
\setcounter{Definition}{0} \setcounter{Remark}{0}

In the present section we formulate the randomized infinite horizon optimal control problem. Firstly, we fix a finite positive measure $\lambda$ on $(A,\Bc(A))$ with full topological support. We also fix a deterministic point $a_0$ in $A$.

Let $(\bar\Omega,\bar\Fc,\bar\P)$ be a complete probability space on which a $d$-dimensional Brownian motion $\bar W=(\bar W_t)_{t\geq0}$ and a Poisson random measure $\bar\mu$ on $[0,\infty)\times A$ are defined. The Poisson random measure $\bar\mu=\sum_{n\geq1}\delta_{(\bar T_n,\bar A_n)}$ is associated with a marked point process $(\bar T_n,\bar A_n)_{n\geq1}$ on $[0,\infty)\times A$, where $(\bar T_n)_{n\geq1}$ is the sequence of jump times, while $(A_n)_{n\geq1}$ is the sequence of $A$-valued marks. The compensator of $\bar\mu$ is given by $\lambda(da)dt$. We denote by $\bar\F^{W,\mu}=(\bar\Fc_t^{W,\mu})_{t\geq0}$ the $\bar\P$-completion of the filtration generated by $\bar W$ and $\bar\mu$, by $\Pc_T(\bar\F^{W,\mu})$, $T\in(0,\infty)$, the predictable $\sigma$-algebra on $[0,T]\times\bar\Omega$ associated with $\bar\F^{W,\mu}$, and by $\Pc(\bar\F^{W,\mu})$ the predictable $\sigma$-algebra on $[0,\infty)\times\bar\Omega$ associated with $\bar\F^{W,\mu}$.

We define the $A$-valued pure-jump process
\begin{equation}\label{I}
\bar I_t \ = \ \sum_{n\ge 0} \bar A_n\,1_{[\bar T_n,\bar T_{n+1})}(t), \qquad \text{for all }t\geq 0,
\end{equation}
with the convention $\bar T_0=0$ and $\bar A_0=a_0$, where $a_0$ is the deterministic point fixed at the beginning of this section. We now consider the following equation:
\begin{equation}\label{sde_random}
\bar X_t \ = \ x_0 + \int_0^t b_s(\bar X,\bar I)\,ds + \int_0^t \sigma_s(\bar X,\bar I)\,d\bar W_s, \qquad \text{for all }t\geq0.
\end{equation}
Under either {\bf (A)} or {\bf (A)'}, there exists a unique $\bar\F^{W,\mu}$-progressive continuous process $\bar X=(\bar X_t)_{t\geq0}$ solution to equation \eqref{sde_random}.

The set $\bar\Vc$ of admissible controls for the randomized problem is given by all $\Pc(\bar\F^{W,\mu})\otimes\Bc(A)$-measurable and bounded maps $\bar\nu\colon\Omega\times \R_+\times A\to (0,\infty)$. We also define, for every $n\in\mathbb N\backslash\{0\}$, the set $\bar\Vc_n:=\{\bar\nu\in\bar\Vc\colon\bar\nu\text{ is bounded by }n\}$. Notice that $\bar\Vc=\cup_n\bar\Vc_n$. For every $\bar\nu\in\Vc$, we consider the corresponding Dol\'eans-Dade exponential process
\begin{align} \nonumber
\kappa_t^{\bar\nu} \ &= \
\Ec_t\left(\int_0^\cdot\int_A (\bar\nu_s(a) - 1)\,(\bar\mu(ds\, da)- \lambda(da)\,ds)\right) \\
&= \ \exp\bigg(\int_0^t\int_A (1 - \bar\nu_s(a))\lambda(da)\,ds
\bigg)\prod_{0<\bar T_n\le t}\bar\nu_{\bar T_n}(\bar A_n),\label{doleans}
\end{align}
for all $t\geq0$. Notice that $\kappa^{\bar\nu}$ is a $(\bar\P,\bar\F^{W,\mu})$-martingale, since $\bar\nu$ is bounded. Then, for every $T>0$, we define the probability
$\bar\P_T^{\bar\nu}(d\bar\omega)=\kappa_T^{\bar\nu}(\bar\omega)\,\bar\P(d\bar\omega)$ on $(\bar\Omega,\bar\Fc_T^{W,\mu})$. For every $T>0$, by
Girsanov's theorem, under $\bar\P_T^{\bar\nu}$
the $(\Fc_t^{W,\mu})_{t\in[0,T]}$-compensator of $\bar\mu$ on $[0,T]\times A$ is $\bar\nu_t(a)\lambda(da)dt$, and $(\bar W_t)_{t\in[0,T]}$ is still a Brownian motion on $[0,T]$ under $\bar\P_T^{\bar\nu}$. Notice that, under {\bf (A)}, for every $T>0$ and $p>0$, we have that there exists a constant $C_{T,p}\geq0$ such that ($\bar\E_T^{\bar\nu}$ denotes the $\bar\P_T^{\bar\nu}$-expectation)
\[
\bar\E_T^{\bar\nu}\Big[\sup_{s\in[0,T]}|\bar X_s|^p\Big] \ \leq \ C_{T,p}\big(1 + |x_0|^p\big).
\]
On the other hand, let {\bf (A)'} hold. Then, by Lemma \ref{L:Estimate} we know that there exists two positive constants $\bar C$ and $\bar\beta$ such that
\[
\bar\E_T^{\bar\nu}\Big[\sup_{s\in[0,T]}|\bar X_s|^r\Big] \ \leq \ \bar Ce^{\bar\beta T}(1 + |x_0|^r),
\]
with $\bar C$ and $\bar\beta$ independent of $T>0$, $\bar\nu\in\bar\Vc$, $x_0\in\R^n$. More generally,  by Lemma \ref{L:Estimate} we have the following estimate:
\begin{equation}\label{beta_rand}
\bar\E_T^{\bar\nu}\Big[\sup_{s\in[0,T]}|\bar X_s|^r\Big|\bar\Fc_t^{W,\mu}\Big] \ \leq \ \bar Ce^{\bar\beta(T-t)}(1 + \sup_{s\in[0,t]}|\bar X_s|^r), \qquad \bar\P\text{-a.s.}
\end{equation}
For all $T>0$, we consider the finite horizon randomized control problem
\[
V_T^{\Rc} \ = \ \sup_{\bar\nu\in\bar\Vc} J_T^\Rc(\bar \nu),
\]
where
\[
J_T^\Rc(\bar\nu) \ = \ \bar\E_T^{\bar\nu}
\bigg[\int_0^T e^{-\beta t}\,f_t(\bar X,\bar I)\,dt\bigg].
\]
Finally, we define the value of the randomized control problem as follows:
\begin{equation}\label{V_random}
V^\Rc \ := \ \lim_{T\rightarrow\infty} V_T^\Rc \ = \ \lim_{T\rightarrow\infty} \sup_{\bar\nu\in\bar\Vc} J_T^\Rc(\bar\nu).
\end{equation}
Proceeding as in the proof of Lemma \ref{L:ApproxV}, it is easy to see that
\[
\sup_{T,T'\geq S} \big|V_T^\Rc - V_{T'}^\Rc\big| \ \overset{S\rightarrow\infty}{\longrightarrow} \ 0,
\]
therefore the limit in \eqref{V_random} exists.
\begin{Remark}\label{girsanovinfinite}
\rm{
Assume that either {\bf (A)} or {\bf (A)'} holds, and suppose that $(\bar\Omega,\bar\Fc,\bar\P)$ has a canonical representation. More precisely, consider the following sets:
\begin{itemize}
\item $\Omega'$ the set of continuous trajectories $\omega'\colon[0,\infty)\rightarrow\R^d$ satisfying $\omega'(0)=0$. We denote $\bar W$ the canonical process on $\Omega'$, $(\Fc_t^W)_{t\geq0}$ the canonical filtration, $\P'$ the Wiener measure on $(\Omega',\Fc_\infty^W)$;
\item $\Omega''$ is the set of double sequences $\omega''=(t_n,a_n)_{n\geq1}\subset(0,\infty)\times A$ satisfying $t_n<t_{n+1}\nearrow\infty$. We denote
 $(\bar T_n,\bar A_n)_{n\geq1}$ the canonical marked point process, $\bar\mu$ $=$ $\sum_{n\ge 1}\delta_{(\bar T_n,\bar A_n)}$ the associated random measure, $(\Fc_t^\mu)_{t\geq0}$ the filtration generated by $\bar\mu$, and $\P''$ the unique probability on
$\Fc_\infty^\mu$ such that $\bar\mu$ is a Poisson random measure with compensator $\lambda(da)dt$.
\end{itemize}
Now, let $\bar\Omega=\Omega'\times\Omega''$, let $\bar \calf$ be the completion of
$\calf_\infty^W\otimes\Fc_\infty^\mu$ with respect to
$\P'\otimes\P''$ and let
$\bar\P$ be the extension of $\P'\otimes\P''$ to $\bar\calf$. Notice that $\bar W$ and $\bar\mu$ can be extended in a canonical way to $\bar\Omega$. We will denote these extensions by
the same symbols. We also denote by $\F^{W,\mu}=(\Fc_t^{W,\mu})_{t\geq0}$ the filtration generated by $\bar W$ and $\bar\mu$, and by $\bar\F^{W,\mu}=(\bar\Fc_t^{W,\mu})_{t\geq0}$ the $\bar\P$-completion of $\F^{W,\mu}$.

Recall that, for every $T>0$, $\bar\P_T^{\bar\nu}$ is a probability on $(\bar\Omega,\bar\Fc_T^{W,\mu})$, then, in particular, on $(\bar\Omega,\Fc_T^{W,\mu})$. Moreover, the following consistency condition holds: $\bar\P_T^{\bar\nu}$ coincides with $\bar\P_t^{\bar\nu}$ on $\Fc_t^{W,\mu}$, whenever $0< t\leq T$. Then, by Kolmogorov's extension theorem, we deduce that there exists a probability measure $\bar\P^{\bar\nu}$ on $(\bar\Omega,\Fc_\infty^{W,\mu})$ such that $\bar\P^{\bar\nu}$ coincides with $\bar\P_T^{\bar\nu}$ on $\Fc_T^{W,\mu}$, for all $T>0$.

Notice that $\bar\P^{\bar\nu}$ can be defined in a consistent way only on $\Fc_\infty^{W,\mu}$ (rather than on $\bar\Fc_\infty^{W,\mu}$). Indeed, since the martingale $(\kappa_t^{\bar\nu})_{t\geq0}$ is in general not uniformly integrable, it follows from Proposition VIII.1.1 in \cite{revyor99} that $\bar\P^{\bar\nu}$ is in general not absolutely continuous with respect to $\bar\P$ on $\Fc_\infty^{W,\mu}$. In particular, when this is the case, there exists some $N\in\Fc_\infty^{W,\mu}$ such that $N\notin\Fc_T^{W,\mu}$, for all $T>0$, and $\bar\P(N)=0$, however $\bar\P^{\bar\nu}(N)>0$. Since $N\in\bar\Fc_T^{W,\mu}$ for all $T>0$, we have $\bar\P_T^{\bar\nu}(N)=0$, therefore we can not extend $\bar\P^{\bar\nu}$ to $\bar\Fc_\infty^{W,\mu}$ without violating the consistency condition: $\bar\P^{\bar\nu}=\bar\P_T^{\bar\nu}$ on $\bar\Fc_T^{W,\mu}$.

Now, notice that the process $\bar I$ is given by \eqref{I} and hence it is $\F^{W,\mu}$-adapted. On the other hand, the process $\bar X$, solution to equation \eqref{sde_random}, is $\bar\F^{W,\mu}$-progressive and continuous, therefore it is $\bar\F^{W,\mu}$-predictable. By IV-78 in \cite{DellMeyerA} it follows that there exists an $\F^{W,\mu}$-predictable process $\hat X$, such that $\bar X$ and $\hat X$ are $\bar\P$-indistinguishable. Then, we have the following representation for $V^\Rc$:
\begin{equation}\label{value_random}
V^\Rc \ = \ \sup_{\bar\nu\in\bar\Vc} J^\Rc(\bar \nu),
\end{equation}
where ($\bar\E^{\bar\nu}$ denotes the $\bar\P^{\bar\nu}$-expectation)
\begin{equation} \label{J_random}
J^\Rc(\bar\nu) \ := \ \bar\E^{\bar\nu}
\bigg[\int_0^\infty e^{-\beta t}\,f_t(\hat X,\bar I)\,dt\bigg].
\end{equation}
Let us prove formula \eqref{value_random}. We begin noting that, for all $T>0$,
\[
J_T^\Rc(\bar\nu) \ = \ \bar\E_T^{\bar\nu}
\bigg[\int_0^T e^{-\beta t}\,f_t(\bar X,\bar I)\,dt\bigg] \ = \ \bar\E_T^{\bar\nu}
\bigg[\int_0^T e^{-\beta t}\,f_t(\hat X,\bar I)\,dt\bigg] \ = \ \bar\E^{\bar\nu}
\bigg[\int_0^T e^{-\beta t}\,f_t(\hat X,\bar I)\,dt\bigg].
\]
Then, under either {\bf (A)} or {\bf (A)'}, proceeding along the same lines as in the proof of Lemma \ref{L:ApproxV}, we obtain
\[
\sup_{\bar\nu\in\bar\Vc} \big|J^\Rc(\bar \nu) - J_T^\Rc(\bar \nu)\big| \ \overset{T\rightarrow\infty}{\longrightarrow} \ 0
\]
and so
\[
\lim_{T\rightarrow\infty} V_T^{\Rc} \ = \ \sup_{\bar\nu\in\bar\Vc} J^\Rc(\bar\nu).
\]
Since, by definition, $V^\Rc=\lim_{T\rightarrow\infty} V_T^{\Rc}$, we conclude that formula \eqref{value_random} holds.
\ep}
\end{Remark}

\section{Identification of the values and backward SDE representation}
\label{S:BSDE}

\setcounter{equation}{0} \setcounter{Assumption}{0}
\setcounter{Theorem}{0} \setcounter{Proposition}{0}
\setcounter{Corollary}{0} \setcounter{Lemma}{0}
\setcounter{Definition}{0} \setcounter{Remark}{0}

In the present section we prove that the original control problem and the randomized control problem have the same value, namely $V=V^\Rc$. We exploit this result in order to derive a backward stochastic differential equation representation for the value $V$.

\begin{Theorem}\label{T:Identification}
Under either {\bf (A)} or {\bf (A')}, we have
\begin{equation}\label{V_T=V_T_random}
V_T \ = \ V_T^\Rc,
\end{equation}
for all $T\in(0,\infty)$, and also
\begin{equation}\label{V=V_random}
V \ = \ V^\Rc.
\end{equation}
\end{Theorem}
\textbf{Proof.}
We begin noting that identity \eqref{V=V_random} is a straightforward consequence of identity \eqref{V_T=V_T_random}, Lemma \ref{L:ApproxV}, and definition \eqref{V_random}. On the other hand, for every $T\in(0,\infty)$, identity \eqref{V_T=V_T_random}
is a direct consequence of Theorem 3.1 in \cite{BCFP16b}.
\ep

\vspace{3mm}

We now prove that $V$ is related to the following infinite horizon backward stochastic differential equation with nonpositive jumps:
\begin{align}
Y_t \ &= \ Y_T - \beta\int_t^T Y_s\,ds + \int_t^T f_s(\bar X,\bar I)\,ds + K_T - K_t \notag \\
&\quad \ - \int_t^T Z_s\,d\bar W_s - \int_t^T \int_A U_s(a)\,\bar\mu(ds,da), \qquad 0 \leq t \leq T, \;  \forall\,T \in (0,\infty), \; \bar\P\text{-a.s.} \label{BSDE} \\
\notag \\
U_t(a) \ &\leq \ 0, \qquad dt\otimes d\bar\P\otimes\lambda(da)\text{-a.e. on }[0,\infty)\times\bar\Omega\times A. \label{JumpConstr}
\end{align}

Let us introduce some additional notations. Given $T\in(0,\infty)$, we denote:
\begin{itemize}
\item ${\bf S^\infty}$, the family of real c\`adl\`ag $\bar\F^{W,\mu}$-adapted stochastic processes $Y$ $=$
$(Y_t)_{t\geq0}$ on $(\bar\Omega,\bar\Fc,\bar\P)$ which are uniformly bounded.
\item ${\bf S^2(0,T)}$, the family of real c\`adl\`ag $\bar\F^{W,\mu}$-adapted stochastic processes $Y$ $=$
$(Y_t)_{0\leq t\leq T}$ on $(\bar\Omega,\bar\Fc,\bar\P)$, with
\[
\|Y\|_{_{{\bf S^2(0,T)}}}^2 \ := \ \bar\E\Big[\sup_{0\leq t\leq T}|Y_t|^2\Big] \ < \ \infty.
\]
We set ${\bf S_{\textup{loc}}^2}:=\cap_{T>0}{\bf S^2(0,T)}$.
\item ${\bf L^2(W;0,T)}$, the family of $\R^d$-valued $\Pc_T(\bar\F^{W,\mu})$-measurable stochastic processes
$Z=(Z_t)_{0\leq t\leq T}$ on $(\bar\Omega,\bar\Fc,\bar\P)$, with
\[
\|Z\|_{_{{\bf L^2(W;0,T)}}}^2 \ := \ \bar\E\bigg[\int_0^T |Z_t|^2\,dt\bigg] \ < \ \infty.
\]
We set ${\bf L_{\textup{loc}}^2(W)}:=\cap_{T>0}{\bf L^2(W;0,T)}$.
\item ${\bf L^2(\bar\mu;0,T)}$, the family of
$\Pc_T(\bar\F^{W,\mu})\otimes\Bc(A)$-measurable maps $U\colon[0,T]\times\bar\Omega\times A\rightarrow \R$ such that
\[
\|U\|_{_{{\bf L^2(\bar\mu;0,T)}}}^2 \ := \ \bar\E\bigg[ \int_0^T\int_A  |U_t(a)|^2 \,\lambda(da)dt\bigg] \ < \ \infty.
\]
We set ${\bf L_{\textup{loc}}^2(\bar\mu)}:=\cap_{T>0}{\bf L^2(\bar\mu;0,T)}$.
\item ${\bf K^2(0,T)}$, the family of  nondecreasing c\`adl\`ag $\Pc_T(\bar\F^{W,\mu})$-measurable stochastic processes $K$ $=$ $(K_t)_{0\leq t\leq T}$ such that $\bar\E[|K_T|^2]$ $<$ $\infty$ and $K_0$ $=$ $0$. We set ${\bf K_{\textup{loc}}^2}:=\cap_{T>0}{\bf K^2(0,T)}$.
\end{itemize}

\begin{Theorem}\label{T:Main}
Under {\bf (A)} $($resp. {\bf (A')}$)$ there exists a unique quadruple $(\bar Y,\bar Z,\bar U,\bar K)$ in ${\bf S_{\textup{loc}}^2}\times{\bf L_{\textup{loc}}^2(W)}\times{\bf L_{\textup{loc}}^2(\bar\mu)}\times{\bf K_{\textup{loc}}^2}$, such that $\bar Y\in{\bf S^\infty}$ $($resp. $|\bar Y_t|\leq C(1+\sup_{s\in[0,t]}|\bar X_s|^r)$, for all $t\geq0$, $\bar\P$-a.s., and for some positive constant $C$$)$, satisfying \eqref{BSDE}-\eqref{JumpConstr} which is minimal in the following sense: for any other quadruple $(\hat Y,\hat Z,\hat U,\hat K)$ in ${\bf S_{\textup{loc}}^2}\times{\bf L_{\textup{loc}}^2(W)}\times{\bf L_{\textup{loc}}^2(\bar\mu)}\times{\bf K_{\textup{loc}}^2}$, with $\hat Y\in{\bf S^\infty}$ $($resp. $|\hat Y_t|\leq C(1+\sup_{s\in[0,t]}|\bar X_s|^r)$, for all $t\geq0$, $\bar\P$-a.s., and for some positive constant $C$$)$, satisfying \eqref{BSDE}-\eqref{JumpConstr} we have:
\[
\bar Y_t \ \leq \ \hat Y_t, \qquad \text{for all }t\geq0, \; \bar\P\text{-a.s.}
\]
Under either {\bf (A)} or {\bf (A')}, we have $V=\bar Y_0$ $\bar\P$-a.s.. Moreover, the following formula holds:
\begin{equation}\label{DPP}
\bar Y_t \ = \ \esssup_{\bar\nu\in\bar\Vc}\bar\E_T^{\bar\nu}\bigg[\int_t^\tau e^{-\beta(s-t)}\,f_s(\bar X,\bar I)\,ds + e^{-\beta(\tau-t)}\,\bar Y_\tau\bigg|\bar\Fc_t^{W,\mu}\bigg], \qquad \bar\P\text{-a.s.}
\end{equation}
for all $T>0$, $t\in[0,T]$, and any $\bar\F^{W,\mu}$-stopping time $\tau$ taking values in $[t,T]$.
\end{Theorem}
\textbf{Proof.} We split the proof into four steps.

\vspace{2mm}

\noindent\textbf{Step I.} \emph{Constrained BSDE on $[0,T]$ and corresponding penalized BSDE.} For every $T\in(0,\infty)$, it follows from Theorem 5.1 in \cite{BCFP16b} that $V_T=Y_0^T$ $\bar\P$-a.s.,
where $(Y^T,Z^T,U^T,K^T)$ in ${\bf S^2(0,T)}\times{\bf L^2(W;0,T)}\times{\bf L^2(\bar\mu;0,T)}\times{\bf K^2(0,T)}$ is the unique quadruple satisfying the following backward stochastic differential equation with nonpositive jumps on $[0,T]$ (with zero terminal condition at time $T$):
\begin{align}
Y_t^T \ &= \ - \beta\int_t^T Y_s^T\,ds + \int_t^T f_s(\bar X,\bar I)\,ds + K_T^T - K_t^T \notag \\
&\quad \ - \int_t^T Z_s^T\,d\bar W_s - \int_t^T \int_A U_s^T(a)\,\bar\mu(ds,da), \qquad 0 \leq t \leq T, \; \bar\P\text{-a.s.} \label{BSDE_T} \\
\notag \\
U_t^T(a) \ &\leq \ 0, \qquad dt\otimes d\bar\P\otimes\lambda(da)\text{-a.e. on }[0,T]\times\bar\Omega\times A, \label{JumpConstr_T}
\end{align}
which is minimal in the following sense: for any other quadruple $(\hat Y^T,\hat Z^T,\hat U^T,\hat K^T)$ in ${\bf S^2(0,T)}\times{\bf L^2(W;0,T)}\times{\bf L^2(\bar\mu;0,T)}\times{\bf K^2(0,T)}$ satisfying \eqref{BSDE_T}-\eqref{JumpConstr_T} we have:
\[
Y_t^T \ \leq \ \hat Y_t^T, \qquad \text{for all }0\leq t\leq T, \; \bar\P\text{-a.s.}
\]
Now, for every $n\in\mathbb N$, consider the following penalized backward stochastic differential equation:
\begin{align}
Y_t^{T,n} \ &= \ - \beta\int_t^T Y_s^{T,n}\,ds + \int_t^T f_s(\bar X,\bar I)\,ds + K_T^{T,n} - K_t^{T,n} \notag \\
&\quad \ - \int_t^T Z_s^{T,n}\,d\bar W_s - \int_t^T \int_A U_s^{T,n}(a)\,\bar\mu(ds,da), \qquad 0 \leq t \leq T, \; \bar\P\text{-a.s.} \label{BSDE_T,n}
\end{align}
where
\[
K_t^{T,n} \ = \ n\,\int_0^t\int_A \big(U_s^{T,n}(a)\big)^+ \, \lambda(da)ds.
\]
It is well-known (see Lemma 2.4 in \cite{tang_li}) that there exists a unique triplet $(Y^{T,n},Z^{T,n},U^{T,n})$ in ${\bf S^2(0,T)}\times{\bf L^2(W;0,T)}\times{\bf L^2(\bar\mu;0,T)}$ satisfying \eqref{BSDE_T,n}.
Moreover, we recall from the proof of Theorem 5.1 in \cite{BCFP16b}, that, for every $0\leq t\leq T$, $Y_t^T$ is the $\bar\P$-a.s. increasing limit of the sequence $(Y_t^{T,n})_{n\in\mathbb N}$.

Our aim is now to pass to the limit in \eqref{BSDE_T,n} as $T\rightarrow\infty$, for any fixed $n\in\N$.

\vspace{2mm}

\noindent\textbf{Substep I.a.} \emph{Convergence of $(Y^{T,n})_{T>0}$.} Applying It\^o's formula to the process $(e^{-\beta t}\,Y_t^{T,n})_{t\in[0,T]}$, we obtain
\begin{align}
e^{-\beta t}\,Y_t^{T,n} \ &= \ \int_t^T e^{-\beta s}\,f_s(\bar X,\bar I)\,ds + \int_t^T e^{-\beta s}\,dK_s^{T,n} \notag \\
&\quad \ - \int_t^T e^{-\beta s}\,Z_s^{T,n}\,d\bar W_s - \int_t^T \int_A e^{-\beta s}\,U_s^{T,n}(a)\,\bar\mu(ds,da). \label{BSDE_T_n}
\end{align}
Given $\bar\nu\in\bar\Vc_n$, notice that
\[
\bar\E_T^{\bar\nu}\bigg[\int_t^T \int_A e^{-\beta s}\,U_s^{T,n}(a)\,\bar\mu(ds,da)\bigg|\bar\Fc_t^{W,\mu}\bigg] \ = \ \bar\E_T^{\bar\nu}\bigg[\int_t^T \int_A e^{-\beta s}\,U_s^{T,n}(a)\,\bar\nu_s(a)\lambda(da)ds\bigg|\bar\Fc_t^{W,\mu}\bigg]
\]
and
\begin{align*}
\bar\E_T^{\bar\nu}\bigg[\int_t^T e^{-\beta s}\,dK_s^{T,n}\bigg|\bar\Fc_t^{W,\mu}\bigg] \ &= \ \bar\E_T^{\bar\nu}\bigg[\int_t^T \int_A n\,e^{-\beta s}\,\big(U_s^{T,n}(a)\big)^+\,\lambda(da)ds\bigg|\bar\Fc_t^{W,\mu}\bigg] \\
&\geq \ \bar\E_T^{\bar\nu}\bigg[\int_t^T \int_A e^{-\beta s}\,U_s^{T,n}(a)\,\bar\nu_s(a)\lambda(da)ds\bigg|\bar\Fc_t^{W,\mu}\bigg],
\end{align*}
where we have used the numerical inequality $n\,u^+\geq\nu\,u$, valid for any real numbers $u$ and $\nu$, with $\nu\in[0,n]$. Then, taking the $\bar\P_T^{\bar\nu}$-conditional expectation with respect to $\bar\Fc_t^{W,\mu}$ in \eqref{BSDE_T_n}, we find
\[
e^{-\beta t}\,Y_t^{T,n} \ \geq \ \bar\E_T^{\bar\nu}\bigg[ \int_t^T e^{-\beta s}\,f_s(\bar X,\bar I)\,ds\bigg|\bar\Fc_t^{W,\mu}\bigg].
\]
Therefore
\begin{equation}\label{dual_Y_T,n_ineq}
Y_t^{T,n} \ \geq \ \esssup_{\bar\nu\in\bar\Vc_n}\bar\E_T^{\bar\nu}\bigg[ \int_t^T e^{-\beta(s-t)}\,f_s(\bar X,\bar I)\,ds\bigg|\bar\Fc_t^{W,\mu}\bigg].
\end{equation}
On the other hand, for every $\eps\in(0,1)$, let $\bar\nu^{T,n,\eps}\in\bar\Vc_n$ be given by
\[
\bar\nu_t^{T,n,\eps}(a) \ = \ n\,1_{\{U_t^{T,n}(a)\geq0\}} + \frac{\eps}{T\lambda(A)}\,1_{\{-1\leq U_t^{T,n}(a)<0\}} + \frac{\eps}{T\lambda(A)U_t^{T,n}(a)}\,1_{\{U_t^{T,n}(a)\leq-1\}}.
\]
Taking the $\bar\P^{\bar\nu^{T,n,\eps}}$-conditional expectation with respect to $\bar\Fc_t^{W,\mu}$ in \eqref{BSDE_T_n}, we obtain
\begin{align*}
e^{-\beta t}\,Y_t^{T,n} \ &\leq \ \bar\E_T^{\bar\nu^{T,n,\eps}}\bigg[\int_t^T e^{-\beta s}\,f_s(\bar X,\bar I)\,ds\bigg|\bar\Fc_t^{W,\mu}\bigg] + \eps \\
&\leq \ \esssup_{\bar\nu\in\bar\Vc_n}\bar\E_T^{\bar\nu}\bigg[\int_t^T e^{-\beta s}\,f_s(\bar X,\bar I)\,ds\bigg|\bar\Fc_t^{W,\mu}\bigg] + \eps.
\end{align*}
By the arbitrariness of $\eps$, and using also inequality \eqref{dual_Y_T,n_ineq}, we conclude that
\begin{equation}\label{dual_Y^T,n}
Y_t^{T,n} \ = \ \esssup_{\bar\nu\in\bar\Vc_n}\bar\E_T^{\bar\nu}\bigg[\int_t^T e^{-\beta(s-t)}\,f_s(\bar X,\bar I)\,ds\bigg|\bar\Fc_t^{W,\mu}\bigg], \qquad \bar\P\text{-a.s., for all }0\leq t\leq T.
\end{equation}
Taking the absolute value of both sides, we obtain,
\begin{equation}\label{Y_Estimate}
\big|Y_t^{T,n}\big| \ \leq \ \esssup_{\bar\nu\in\bar\Vc_n}\bar\E_T^{\bar\nu}\bigg[\int_t^T e^{-\beta(s-t)}\,\big|f_s(\bar X,\bar I)\big|\,ds\bigg|\bar\Fc_t^{W,\mu}\bigg], \qquad \bar\P\text{-a.s., for all }0\leq t\leq T.
\end{equation}
Moreover, for any $T'>0$, from \eqref{dual_Y^T,n} we find
\begin{equation}\label{Y_T',n-Y_T,n}
\big|Y_t^{T',n} - Y_t^{T,n}\big| \leq \esssup_{\bar\nu\in\bar\Vc_n}\bar\E_{T\vee T'}^{\bar\nu}\bigg[\int_{T\wedge T'}^{T\vee T'} e^{-\beta(s-t)}\,\big|f_s(\bar X,\bar I)\big|\,ds\bigg|\bar\Fc_t^{W,\mu}\bigg], \quad \bar\P\text{-a.s., for all }0\leq t\leq T\wedge T'.
\end{equation}
Now, we distinguish two cases.

\begin{itemize}
\item \emph{Assumption {\bf (A)} holds}. Since $f$ is bounded, from \eqref{Y_Estimate} we find
\[
\big|Y_t^{T,n}\big| \ \leq \ \|f\|_\infty\frac{1 - e^{-\beta(T-t)}}{\beta} \ \leq \ \frac{\|f\|_\infty}{\beta},
\]
$\bar\P$-a.s., for all $0\leq t\leq T$. Since $(Y_t^{T,n})_{t\in[0,T]}$ is a c\`adl\`ag process, we obtain
\begin{equation}\label{bound_Y^T,n}
\big|Y_t^{T,n}\big| \ \leq \ \frac{\|f\|_\infty}{\beta},
\end{equation}
for all $0\leq t\leq T$, $\bar\P$-a.s., so that $Y^{T,n}$ is a uniformly bounded process. Proceeding in a similar way, we can deduce from estimate \eqref{Y_T',n-Y_T,n} that
\begin{equation}\label{Conv_Y_bdd2}
\big|Y_t^{T',n} - Y_t^{T,n}\big| \ \leq \ \frac{\|f\|_\infty}{\beta} e^{-\beta(T\wedge T'-t)}, \qquad \text{for all }0\leq t\leq T\wedge T',\;\bar\P\text{-a.s.}
\end{equation}
In particular, for any $S\in(0,T\wedge T')$, we have
\begin{equation}\label{Conv_Y_bdd}
\sup_{t\in[0,S]}\big|Y_t^{T',n} - Y_t^{T,n}\big| \ \leq \ \frac{\|f\|_\infty}{\beta} e^{-\beta(T\wedge T'-S)}, \qquad \bar\P\text{-a.s.}
\end{equation}
\item \emph{Assumption {\bf (A')} holds}. By \eqref{Y_Estimate} we have
\begin{align*}
\big|Y_t^{T,n}\big| \ &\leq \ M\,\esssup_{\bar\nu\in\bar\Vc_n}\bar\E_T^{\bar\nu}\bigg[\int_t^T e^{-\beta(s-t)}\,\Big(1 + \sup_{u\in[0,s]}|\bar X_u|^r\Big)\,ds\bigg|\bar\Fc_t^{W,\mu}\bigg] \\
&\leq \ M\,\int_t^T e^{-\beta(s-t)}\,\Big(1 + \esssup_{\bar\nu\in\bar\Vc_n}\bar\E_T^{\bar\nu}\Big[\sup_{u\in[0,s]}|\bar X_u|^r\Big|\bar\Fc_t^{W,\mu}\Big]\Big)\,ds,
\end{align*}
$\bar\P$-a.s., for all $0\leq t\leq T$. Using estimate \eqref{beta_rand}, we find
\begin{align*}
\big|Y_t^{T,n}\big| \ &\leq \ M\,\int_t^T e^{-\beta(s-t)}\,\Big(1 + \bar Ce^{\bar\beta(s-t)}(1 + \sup_{u\in[0,t]}|\bar X_u|^r)\Big)\,ds \\
&\leq 2\,M\,(1 + \bar C)\,\Big(1 + \sup_{u\in[0,t]}|\bar X_u|^r\Big)\,\int_t^T e^{-(\beta-\bar\beta)(s-t)}\,ds \\
&\leq \ \frac{2\,M\,(1 + \bar C)}{\beta-\bar\beta}\Big(1 + \sup_{u\in[0,t]}|\bar X_u|^r\Big),
\end{align*}
$\bar\P$-a.s., for all $0\leq t\leq T$. Since $(Y_t^{T,n})_{t\in[0,T]}$ is a c\`adl\`ag process, we obtain
\begin{equation}\label{bound_Y^T,n_pol}
\big|Y_t^{T,n}\big| \ \leq \ \frac{2\,M\,(1 + \bar C)}{\beta-\bar\beta}\Big(1 + \sup_{s\in[0,t]}|\bar X_s|^r\Big),
\end{equation}
for all $0\leq t\leq T$, $\bar\P$-a.s.. In a similar way, starting from estimate \eqref{Y_T',n-Y_T,n} we can prove that
\begin{equation}\label{Estimate_Y-Y'}
\big|Y_t^{T',n} - Y_t^{T,n}\big| \ \leq \ \frac{2\,M\,(1 + \bar C)}{\beta-\bar\beta}\Big(1 + \sup_{s\in[0,t]}|\bar X_s|^r\Big)e^{-(\beta-\bar\beta)(T\wedge T'-t)},
\end{equation}
for all $0\leq t\leq T\wedge T'$, $\bar\P$-a.s.. As a consequence, for any $S\in(0,T\wedge T')$,
\begin{equation}\label{Conv_Y_pol}
\sup_{t\in[0,S]}\big|Y_t^{T',n} - Y_t^{T,n}\big| \ \leq \ \frac{2\,M\,(1 + \bar C)}{\beta-\bar\beta}\Big(1 + \sup_{t\in[0,S]}|\bar X_t|^r\Big)e^{-(\beta-\bar\beta)(T\wedge T'-S)}, \qquad \bar\P\text{-a.s.}
\end{equation}

\end{itemize}
Then, under either {\bf (A)} or {\bf (A')}, we see that for every $n\in\N$ there exists a c\`adl\`ag process $(Y_t^n)_{t\geq0}\in{\bf S_{\textup{loc}}^2}$ such that, for any $S>0$, $Y^n$ is the $\bar\P$-a.s. uniform limit on $[0,S]$ as $T\rightarrow\infty$ of the sequence of c\`adl\`ag processes $(Y^{T,n})_{T>S}$. Moreover, under {\bf (A)}, we have from \eqref{bound_Y^T,n}
\begin{equation}\label{bound_Y^n}
\big|Y_t^n\big| \ \leq \ \frac{\|f\|_\infty}{\beta}, \qquad \text{for all }t\geq0,\;\bar\P\text{-a.s.}
\end{equation}
In particular, $Y^n\in{\bf S^\infty}$. On the other hand, under {\bf (A')}, we deduce from \eqref{bound_Y^T,n_pol}
\begin{equation}\label{bound_Y^n_pol}
\big|Y_t^n\big| \ \leq \ \frac{2\,M\,(1 + \bar C)}{\beta - \bar\beta}\Big(1 + \sup_{s\in[0,t]}|\bar X_s|^r\Big), \qquad \text{for all }t\geq0,\;\bar\P\text{-a.s.}
\end{equation}
Furthermore, under either {\bf (A)} or {\bf (A')}, from formula \eqref{dual_Y^T,n} we have
\begin{equation}\label{dual_Y^n}
Y_t^n \ = \ \lim_{T\rightarrow\infty}\esssup_{\bar\nu\in\bar\Vc_n}\bar\E_T^{\bar\nu}\bigg[\int_t^T e^{-\beta(s-t)}\,f_s(\bar X,\bar I)\,ds\bigg|\bar\Fc_t^{W,\mu}\bigg], \qquad \bar\P\text{-a.s., for all }t\geq0.
\end{equation}

\vspace{2mm}

\noindent\textbf{Substep I.b.} \emph{Convergence of $(Z^{T,n},U^{T,n})_{T>0}$.} Take $T'>0$ and $S\in(0,T\wedge T')$. An application of It\^o's formula to $(e^{-\beta t}|Y_t^{T',n}-Y_t^{T,n}|)^2$ between $0$ and $S$, yields (taking also the $\bar\P$-expectation)
\begin{align*}
&\bar\E\bigg[\int_0^S e^{-2\beta s}\,\big|Z_s^{T',n} - Z_s^{T,n}\big|^2\,ds \bigg] + \bar\E\bigg[\int_0^S \int_A e^{-2\beta s}\,\big|U_s^{T',n}(a) - U_s^{T,n}(a)\big|^2\,\lambda(da)ds \bigg] \\
&= \ \bar\E\big[e^{-2\beta S}\,\big|Y_S^{T',n} - Y_S^{T,n}\big|^2\big] - \big|Y_0^{T',n} - Y_0^{T,n}\big|^2 \\
&\quad \ + 2n\,\bar\E\bigg[\int_0^S\int_A e^{-2\beta s}\big(Y_s^{T',n} - Y_s^{T,n}\big)\big(\big(U_s^{T',n}(a)\big)^+ - \big(U_s^{T,n}(a)\big)^+\big)\,\lambda(da)ds\bigg] \\
&\quad \ - 2\,\bar\E\bigg[\int_0^S\int_A e^{-2\beta s}\big(Y_s^{T',n} - Y_s^{T,n}\big)\big(U_s^{T',n}(a) - U_s^{T,n}(a)\big)\,\lambda(da)ds\bigg].
\end{align*}
Then
\begin{align*}
&\bar\E\bigg[\int_0^S e^{-2\beta s}\,\big|Z_s^{T',n} - Z_s^{T,n}\big|^2\,ds \bigg] + \bar\E\bigg[\int_0^S \int_A e^{-2\beta s}\,\big|U_s^{T',n}(a) - U_s^{T,n}(a)\big|^2\,\lambda(da)ds \bigg] \\
&\leq \bar\E\big[e^{-2\beta S}\,\big|Y_S^{T',n} - Y_S^{T,n}\big|^2\big] + 2(n+1)\,\bar\E\bigg[\int_0^S\int_A e^{-2\beta s}\big|Y_s^{T',n} - Y_s^{T,n}\big|\big|U_s^{T',n}(a) - U_s^{T,n}(a)\big|\,\lambda(da)ds\bigg] \\
&\leq \bar\E\big[e^{-2\beta S}\,\big|Y_S^{T',n} - Y_S^{T,n}\big|^2\big] + 2(n+1)^2 \lambda(A)\bar\E\bigg[\int_0^S e^{-2\beta s}\,\big|Y_s^{T',n} - Y_s^{T,n}\big|^2\,ds\bigg] \\
&\quad + \frac{1}{2} \bar\E\bigg[\int_0^S \int_A e^{-2\beta s}\,\big|U_s^{T',n}(a) - U_s^{T,n}(a)\big|^2\,\lambda(da)ds \bigg].
\end{align*}
Therefore
\begin{align*}
&\frac{1}{2}\,e^{-2\beta S}\,\bar\E\bigg[\int_0^S \big|Z_s^{T',n} - Z_s^{T,n}\big|^2\,ds \bigg] + \frac{1}{2}\,e^{-2\beta S}\,\bar\E\bigg[\int_0^S \int_A \big|U_s^{T',n}(a) - U_s^{T,n}(a)\big|^2\,\lambda(da)ds \bigg] \\
&\leq \ \bar\E\bigg[\int_0^S e^{-2\beta s}\,\big|Z_s^{T',n} - Z_s^{T,n}\big|^2\,ds \bigg] + \frac{1}{2}\bar\E\bigg[\int_0^S \int_A e^{-2\beta s}\,\big|U_s^{T',n}(a) - U_s^{T,n}(a)\big|^2\,\lambda(da)ds \bigg] \\
&\leq \ \bar\E\big[e^{-2\beta S}\,\big|Y_S^{T',n} - Y_S^{T,n}\big|^2\big] + 2(n+1)^2 \lambda(A)\bar\E\bigg[\int_0^S e^{-2\beta s}\,\big|Y_s^{T',n} - Y_s^{T,n}\big|^2\,ds\bigg].
\end{align*}
In conclusion, we find
\begin{align*}
&\bar\E\bigg[\int_0^S \big|Z_s^{T',n} - Z_s^{T,n}\big|^2\,ds \bigg] + \bar\E\bigg[\int_0^S \int_A \big|U_s^{T',n}(a) - U_s^{T,n}(a)\big|^2\,\lambda(da)ds \bigg] \\
&\leq \ 2\,\bar\E\big[\big|Y_S^{T',n} - Y_S^{T,n}\big|^2\big] + 4(n+1)^2 \lambda(A)\bar\E\bigg[\int_0^S e^{2\beta(S-s)}\,\big|Y_s^{T',n} - Y_s^{T,n}\big|^2\,ds\bigg].
\end{align*}
Using either \eqref{Conv_Y_bdd2} under {\bf (A)} or \eqref{Estimate_Y-Y'} under {\bf (A')}, we deduce that
\[
\bar\E\bigg[\int_0^S \big|Z_s^{T',n} - Z_s^{T,n}\big|^2\,ds \bigg] + \bar\E\bigg[\int_0^S \int_A \big|U_s^{T',n}(a) - U_s^{T,n}(a)\big|^2\,\lambda(da)ds \bigg] \ \overset{T,T'\rightarrow\infty}{\longrightarrow} \ 0.
\]
In other words, for any $S>0$, the sequence $(Z^{T,n},U^{T,n})_{T>S}$ is a Cauchy sequence in the Hilbert space ${\bf L^2(W;0,S)}\times{\bf L^2(\bar\mu;0,S)}$. It follows that there exists $(Z^n,U^n)\in{\bf L_{\textup{loc}}^2(W)}\times{\bf L_{\textup{loc}}^2(\bar\mu)}$ such that
\begin{equation}\label{Conv_Z,U}
\bar\E\bigg[\int_0^S \big|Z_s^{T,n} - Z_s^n\big|^2\,ds \bigg] + \bar\E\bigg[\int_0^S \int_A \big|U_s^{T,n}(a) - U_s^n(a)\big|^2\,\lambda(da)ds \bigg] \ \overset{T\rightarrow\infty}{\longrightarrow} \ 0.
\end{equation}
Now, take $S\in(0,T)$ and consider equation \eqref{BSDE_T,n} between $t\in[0,S]$ and $S$:
\begin{align*}
Y_t^{T,n} \ &= \ Y_S^{T,n} - \beta\int_t^S Y_s^{T,n}\,ds + \int_t^S f_s(\bar X,\bar I)\,ds + n\int_t^S \int_A \big(U_s^{T,n}(a)\big)^+\,\lambda(da)ds \\
&\quad \ - \int_t^S Z_s^{T,n}\,d\bar W_s - \int_t^S \int_A U_s^{T,n}(a)\,\bar\mu(ds,da).
\end{align*}
Letting $T\rightarrow\infty$, using either \eqref{Conv_Y_bdd} under {\bf (A)} or \eqref{Conv_Y_pol} under {\bf (A')}, and also \eqref{Conv_Z,U}, we obtain
\begin{align}
Y_t^n \ &= \ Y_S^n - \beta\int_t^S Y_s^n\,ds + \int_t^S f_s(\bar X,\bar I)\,ds + n\int_t^S \int_A \big(U_s^n(a)\big)^+\,\lambda(da)ds \notag \\
&\quad \ - \int_t^S Z_s^n\,d\bar W_s - \int_t^S \int_A U_s^n(a)\,\bar\mu(ds,da), \qquad\qquad 0\leq t\leq S,\;\bar\P\text{-a.s.} \label{BSDE_n}
\end{align}
Since $S$ is arbitrary in \eqref{BSDE_n}, we conclude that $(Y^n,Z^n,U^n)$ is a solution to the above infinite horizon backward stochastic differential equation.

Our aim is now to pass to the limit in \eqref{BSDE_n} as $n\rightarrow\infty$.

\vspace{2mm}

\noindent\textbf{Step II.} \emph{Proof of formula \eqref{DPP}.} Recalling that $\bar\Vc_n\subset\bar\Vc_{n+1}$, by formula \eqref{dual_Y^n} we see that $(Y_t^n)_{n\in\N}$ is an increasing sequence, for all $t\geq0$. In particular, we have $Y_t^0 \leq Y_t^1\leq\cdots\leq Y_t^n\leq\cdots$, $\bar\P$-a.s., for all $t\geq0$. Since $Y^n$, for every $n\in\N$, is a c\`adl\`ag process, we deduce that $Y_t^0 \leq Y_t^1\leq\cdots\leq Y_t^n\leq\cdots$, for all $t\geq0$, $\bar\P$-a.s.. Therefore, there exists an $\bar\F^{W,\mu}$-adapted process $(\bar Y_t)_{t\geq0}$ such that $Y_t^n$ converges pointwise increasingly to $\bar Y_t$, for all $t\geq0$, $\bar\P$-a.s.. Moreover, under {\bf (A)} we have, using estimate \eqref{bound_Y^n},
\begin{equation}\label{bound_Y}
|\bar Y_t| \ \leq \ \frac{\|f\|_\infty}{\beta}, \qquad \text{for all }t\geq0,\;\bar\P\text{-a.s.}
\end{equation}
Therefore $\bar Y$ is uniformly bounded. On the other hand, under {\bf (A')} we obtain, using estimate \eqref{bound_Y^n_pol},
\begin{equation}\label{bound_Y_pol}
|\bar Y_t| \ \leq \ \frac{2\,M\,(1 + \bar C)}{\beta - \bar\beta}\Big(1 + \sup_{s\in[0,t]}|\bar X_s|^r\Big), \qquad \text{for all }t\geq0,\;\bar\P\text{-a.s.}
\end{equation}
Let us now prove formula \eqref{DPP}. Fix $T>0$, $t\in[0,T]$, and an $\bar\F^{W,\mu}$-stopping time $\tau$ taking values in $[t,T]$. We begin noting that, considering equation \eqref{BSDE_n} written between $t$ and $\tau$, and proceeding along the same lines as in the proof of formula \eqref{dual_Y^T,n}, taking into account that the terminal condition is now given by $Y_\tau^n$, we can prove that
\begin{equation}\label{dual_Y^n_bis}
Y_t^n \ = \ \esssup_{\bar\nu\in\bar\Vc_n}\bar\E_T^{\bar\nu}\bigg[\int_t^\tau e^{-\beta(s-t)}\,f_s(\bar X,\bar I)\,ds + e^{-\beta(\tau-t)}\,Y_\tau^n\bigg|\bar\Fc_t^{W,\mu}\bigg], \qquad \bar\P\text{-a.s.}
\end{equation}
Since $\bar\Vc_n\subset\bar\Vc$ and $Y_\tau^n\leq\bar Y_\tau$, $\bar\P$-a.s., we get
\[
Y_t^n \ \leq \ \esssup_{\bar\nu\in\bar\Vc}\bar\E_T^{\bar\nu}\bigg[\int_t^\tau e^{-\beta(s-t)}\,f_s(\bar X,\bar I)\,ds + e^{-\beta(\tau-t)}\,\bar Y_\tau\bigg|\bar\Fc_t^{W,\mu}\bigg], \qquad \bar\P\text{-a.s.}
\]
Recalling that $Y_t^n\nearrow\bar Y_t$, $\bar\P$-a.s., we obtain the inequality
\begin{equation}\label{DPP_Ineq}
\bar Y_t \ \leq \ \esssup_{\bar\nu\in\bar\Vc}\bar\E_T^{\bar\nu}\bigg[\int_t^\tau e^{-\beta(s-t)}\,f_s(\bar X,\bar I)\,ds + e^{-\beta(\tau-t)}\,\bar Y_\tau\bigg|\bar\Fc_t^{W,\mu}\bigg], \qquad \bar\P\text{-a.s.}
\end{equation}
On the other hand, let $n,m\in\N$, with $n\geq m$, then
\begin{align*}
\bar Y_t \ \geq \ Y_t^n \ &= \ \esssup_{\bar\nu\in\bar\Vc_n}\bar\E_T^{\bar\nu}\bigg[\int_t^\tau e^{-\beta(s-t)}\,f_s(\bar X,\bar I)\,ds + e^{-\beta(\tau-t)}\,Y_\tau^n\bigg|\bar\Fc_t^{W,\mu}\bigg] \\
&\geq \ \esssup_{\bar\nu\in\bar\Vc_n}\bar\E_T^{\bar\nu}\bigg[\int_t^\tau e^{-\beta(s-t)}\,f_s(\bar X,\bar I)\,ds + e^{-\beta(\tau-t)}\,Y_\tau^m\bigg|\bar\Fc_t^{W,\mu}\bigg], \qquad \bar\P\text{-a.s.}
\end{align*}
Taking the supremum over $n\in\{m,m+1,\ldots\}$, we find
\[
\bar Y_t \ \geq \ \esssup_{\bar\nu\in\bar\Vc}\bar\E_T^{\bar\nu}\bigg[\int_t^\tau e^{-\beta(s-t)}\,f_s(\bar X,\bar I)\,ds + e^{-\beta(\tau-t)}\,Y_\tau^m\bigg|\bar\Fc_t^{W,\mu}\bigg], \qquad \bar\P\text{-a.s.}
\]
In particular, we have
\[
\bar Y_t \ \geq \ \bar\E_T^{\bar\nu}\bigg[\int_t^\tau e^{-\beta(s-t)}\,f_s(\bar X,\bar I)\,ds + e^{-\beta(\tau-t)}\,Y_\tau^m\bigg|\bar\Fc_t^{W,\mu}\bigg], \qquad \bar\P\text{-a.s.}
\]
for all $\bar\nu\in\bar\Vc$, $m\in\N$. Taking the limit as $m\rightarrow\infty$, and afterwards the $\esssup_{\bar\nu\in\bar\Vc}$, we conclude that
\[
\bar Y_t \ \geq \ \esssup_{\bar\nu\in\bar\Vc}\bar\E_T^{\bar\nu}\bigg[\int_t^\tau e^{-\beta(s-t)}\,f_s(\bar X,\bar I)\,ds + e^{-\beta(\tau-t)}\,\bar Y_\tau\bigg|\bar\Fc_t^{W,\mu}\bigg], \qquad \bar\P\text{-a.s.}
\]
which, together with \eqref{DPP_Ineq}, gives formula \eqref{DPP}.

\vspace{2mm}

\noindent\textbf{Step III.} \emph{Convergence of the penalized infinite horizon BSDE.} By either \eqref{bound_Y} or \eqref{bound_Y_pol}, we see that
\[
\esssup_{\bar\nu\in\bar\Vc}\bar\E_T^{\bar\nu}\big[e^{-\beta(T-t)}|\bar Y_T|\big|\bar\Fc_t^{W,\mu}\big] \ \overset{T\rightarrow\infty}{\underset{\bar\P\text{-a.s.}}{\longrightarrow}} \ 0.
\]
Then, from \eqref{DPP} with $\tau=T$, we obtain
\begin{equation}\label{dual_Y}
\bar Y_t \ = \ \lim_{T\rightarrow\infty}\esssup_{\bar\nu\in\bar\Vc}\bar\E_T^{\bar\nu}\bigg[\int_t^T e^{-\beta(s-t)}\,f_s(\bar X,\bar I)\,ds\bigg|\bar\Fc_t^{W,\mu}\bigg], \qquad \bar\P\text{-a.s.}
\end{equation}
In particular, taking $t=0$ in \eqref{dual_Y}, we see that $V=\bar Y_0$, $\bar\P$-a.s..
Now, let $T>0$ and consider the following backward stochastic differential equation with nonpositive jumps on $[0,T]$ with terminal condition $\bar Y_T$:
\begin{align}
Y_t^T \ &= \ \bar Y_T - \beta\int_t^T Y_s^T\,ds + \int_t^T f_s(\bar X,\bar I)\,ds + K_T^T - K_t^T \notag \\
&\quad \ - \int_t^T Z_s^T\,d\bar W_s - \int_t^T \int_A U_s^T(a)\,\bar\mu(ds,da), \qquad\qquad 0\leq t\leq T,\;\bar\P\text{-a.s.} \label{BSDE_T_bis} \\
\notag \\
U_t^T(a) \ &\leq \ 0, \qquad dt\otimes d\bar\P\otimes\lambda(da)\text{-a.e. on }[0,T]\times\bar\Omega\times A. \label{Jump_T_bis}
\end{align}
By Theorem 2.1 in \cite{KP12}, we know that there exists a unique minimal solution $(Y^T,Z^T,U^T,K^T)\in{\bf S^2(0,T)}\times{\bf L^2(W;0,T)}\times{\bf L^2(\bar\mu;0,T)}\times{\bf K^2(0,T)}$ to equation \eqref{BSDE_T_bis}-\eqref{Jump_T_bis}, where minimal means that given any other solution $(\hat Y^T,\hat Z^T,\hat U^T,\hat K^T)\in{\bf S^2(0,T)}\times{\bf L^2(W;0,T)}\times{\bf L^2(\bar\mu;0,T)}\times{\bf K^2(0,T)}$ to equation \eqref{BSDE_T_bis}-\eqref{Jump_T_bis}, then
\[
Y_t^T \ \leq \ \hat Y_t^T, \qquad \text{for all }0\leq t\leq T, \; \bar\P\text{-a.s.}
\]
On the other hand, consider the penalized infinite horizon backward stochastic differential equation \eqref{BSDE_n} as an equation on $[0,T]$:
\begin{align*}
Y_t^n \ &= \ Y_T^n - \beta\int_t^T Y_s^n\,ds + \int_t^T f_s(\bar X,\bar I)\,ds + n\int_t^T \int_A \big(U_s^n(a)\big)^+\,\lambda(da)ds \\
&\quad \ - \int_t^T Z_s^n\,d\bar W_s - \int_t^T \int_A U_s^n(a)\,\bar\mu(ds,da), \qquad\qquad 0\leq t\leq T,\;\bar\P\text{-a.s.}
\end{align*}
Then, proceeding as in the proof of Theorem 2.1 in \cite{KP12}, we can prove that:
\begin{itemize}
\item[(i)] for all $0\leq t\leq T$, $Y_t^T$ is the increasing pointwise $\bar\P$-a.s. limit of the sequence $(Y_t^n)_{n\in\N}$;
\item[(ii)] $(Z^T,U^T)$ is the weak limit in ${\bf L^2(W;0,T)}\times{\bf L^2(\bar\mu;0,T)}$ of the sequence $(Z_{|[0,T]}^n,U_{|[0,T]}^n)_{n\in\N}$;
\item[(iii)] for all $0\leq t\leq T$, $K_t^T$ is the weak limit in $L^2(\Omega,\bar\Fc_t^{W,\mu},\bar\P)$ of the sequence $(K_t^n)_{n\in\N}$.
\end{itemize}
Notice that the only difference with respect to the case considered in \cite{KP12} is that here the penalized equation has a terminal condition $Y_T^n$ depending on $n$; however, the proof in \cite{KP12} still works using the property that $Y_T^n$ converges pointwise increasingly $\bar\P$-a.s. to $\bar Y_T$.

From point (i), we deduce that $\bar Y_t= Y_t^T$, $\bar\P$-a.s., for all $0\leq t\leq T$, and any $T>0$. In particular, we see that $\bar Y$ admits a c\`adl\`ag version. From now on, by an abuse of notation, we denote by $\bar Y$ this c\`adl\`ag version. Then, $\bar Y\in{\bf S_{\textup{loc}}^2}$, and in addition $\bar Y\in{\bf S^\infty}$ under {\bf (A)}. We also see that there exists $(\bar Z,\bar U,\bar K)\in{\bf L_{\textup{loc}}^2(W)}\times{\bf L_{\textup{loc}}^2(\bar\mu)}\times{\bf K_{\textup{loc}}^2}$ such that $\|\bar Z_{|[0,T]}-Z^T\|_{_{{\bf L^2(W;0,T)}}}=0$, $\|\bar U_{|[0,T]}-U^T\|_{_{{\bf L^2(\bar\mu;0,T)}}}=0$, $\|\bar K_{|[0,T]}-K^T\|_{_{{\bf S^2(0,T)}}}=0$. Since, for any $T>0$, $(Y^T,Z^T,U^T,K^T)$ satisfies equation \eqref{BSDE_T_bis}-\eqref{Jump_T_bis}, it follows that $(\bar Y,\bar Z,\bar U,\bar K)$ solves \eqref{BSDE}-\eqref{JumpConstr}.

Finally, let us prove that $(\bar Y,\bar Z,\bar U,\bar K)$ is a minimal solution to equation \eqref{BSDE}-\eqref{JumpConstr}. Consider an arbitrary quadruple $(\hat Y,\hat Z,\hat U,\hat K)\in{\bf S_{\textup{loc}}^2}\times{\bf L_{\textup{loc}}^2(W)}\times{\bf L_{\textup{loc}}^2(\bar\mu)}\times{\bf K_{\textup{loc}}^2}$ solution to equation \eqref{BSDE}-\eqref{JumpConstr}, with either $\hat Y\in{\bf S^\infty}$ under {\bf (A)} or $|\hat Y_t|\leq C(1+\sup_{s\in[0,t]}|\bar X_s|^r)$, for all $t\geq0$, $\bar\P$-a.s., and for some positive constant $C$, under {\bf (A')}. In particular, we have, for any $0\leq t\leq T$,
\begin{align*}
&\hat Y_t + \int_t^T e^{-\beta(s-t)}\,\hat Z_s\,d\bar W_s + \int_t^T \int_A e^{-\beta(s-t)}\,\hat U_s(a)\,\bar\mu(ds,da) \\
&= \ e^{-\beta(T-t)}\,\hat Y_T + \int_t^T e^{-\beta(s-t)}\,f_s(\bar X,\bar I)\,ds + \int_t^T e^{-\beta(s-t)}\,d\hat K_s.
\end{align*}
Given $\bar\nu\in\bar\Vc$, taking the $\bar\P_T^{\bar\nu}$-conditional expectation with respect to $\bar\Fc_t^{W,\mu}$, we obtain
\begin{align*}
\hat Y_t \ &\geq \ \hat Y_t + \bar\E_T^{\bar\nu}\bigg[\int_t^T \int_A e^{-\beta(s-t)}\,\hat U_s(a)\,\lambda(da)ds\bigg|\bar\Fc_t^{W,\mu}\bigg] \\
&= \ \bar\E_T^{\bar\nu}\big[e^{-\beta(T-t)}\,\hat Y_T\big|\bar\Fc_t^{W,\mu}\big] + \bar\E_T^{\bar\nu}\bigg[\int_t^T e^{-\beta(s-t)}\,f_s(\bar X,\bar I)\,ds\bigg|\bar\Fc_t^{W,\mu}\bigg] + \bar\E_T^{\bar\nu}\bigg[\int_t^T e^{-\beta(s-t)}\,d\hat K_s\bigg|\bar\Fc_t^{W,\mu}\bigg] \\
&\geq \ \bar\E_T^{\bar\nu}\big[e^{-\beta(T-t)}\,\hat Y_T\big|\bar\Fc_t^{W,\mu}\big] + \bar\E_T^{\bar\nu}\bigg[\int_t^T e^{-\beta(s-t)}\,f_s(\bar X,\bar I)\,ds\bigg|\bar\Fc_t^{W,\mu}\bigg].
\end{align*}
From the arbitrariness of $\bar\nu$, we get
\[
\hat Y_t \ \geq \ \esssup_{\bar\nu\in\bar\Vc} \bar\E_T^{\bar\nu}\bigg[\int_t^T e^{-\beta(s-t)}\,f_s(\bar X,\bar I)\,ds + e^{-\beta(T-t)}\,\hat Y_T\bigg|\bar\Fc_t^{W,\mu}\bigg].
\]
Using the bounds satisfied by $\hat Y$ under {\bf (A)} or {\bf (A')}, and estimate \eqref{beta_rand}, we see that
\[
\esssup_{\bar\nu\in\bar\Vc}\bar\E_T^{\bar\nu}\big[e^{-\beta(T-t)}|\hat Y_T|\big|\bar\Fc_t^{W,\mu}\big] \ \overset{T\rightarrow\infty}{\underset{\bar\P\text{-a.s.}}{\longrightarrow}} \ 0.
\]
This implies that
\begin{align*}
\hat Y_t \ \geq \ \liminf_{T\rightarrow\infty}\esssup_{\bar\nu\in\bar\Vc}\bar\E_T^{\bar\nu}\bigg[\int_t^T e^{-\beta(s-t)}\,f_s(\bar X,\bar I)\,ds + e^{-\beta(T-t)}\hat Y_T\bigg|\bar\Fc_t^{W,\mu}\bigg]& \\
= \ \liminf_{T\rightarrow\infty}\esssup_{\bar\nu\in\bar\Vc}\bar\E_T^{\bar\nu}\bigg[\int_t^T e^{-\beta(s-t)}\,f_s(\bar X,\bar I)\,ds\bigg|\bar\Fc_t^{W,\mu}\bigg]& \\
\geq \ \lim_{T\rightarrow\infty}\esssup_{\bar\nu\in\bar\Vc_n}\bar\E_T^{\bar\nu}\bigg[\int_t^T e^{-\beta(s-t)}\,f_s(\bar X,\bar I)\,ds\bigg|\bar\Fc_t^{W,\mu}\bigg]& \ = \ Y_t^n, \qquad \bar\P\text{-a.s.}
\end{align*}
Letting $n\rightarrow\infty$, we end up with $\hat Y_t\geq\bar Y_t$, $\bar\P$-a.s., for all $t\geq0$, which yields the minimality of $\bar Y$.

\vspace{2mm}

\noindent\textbf{Step IV.} \emph{Uniqueness.} Since $(\bar Y,\bar Z,\bar U,\bar K)$ is a minimal solution to equation \eqref{BSDE}-\eqref{JumpConstr}, by definition, the component $\bar Y$ is uniquely determined. In order to prove the uniqueness of the other components $(\bar Z,\bar U,\bar K)$, we can proceed as in Remark 2.1 in \cite{KP12} and prove that their restrictions $(\bar Z_{|[0,T]},\bar U_{|[0,T]},\bar K_{|[0,T]})$ to the interval $[0,T]$, for any $T>0$, are uniquely determined. Then, the result follows from the arbitrariness of $T$.
\ep

\section{Feynman-Kac representation for fully non-linear elliptic PDE}
\label{S:Markov}

\setcounter{equation}{0} \setcounter{Assumption}{0}
\setcounter{Theorem}{0} \setcounter{Proposition}{0}
\setcounter{Corollary}{0} \setcounter{Lemma}{0}
\setcounter{Definition}{0} \setcounter{Remark}{0}

We now apply the results of the previous sections in order to determine a non-linear Feynman-Kac representation formula for a fully non-linear elliptic partial differential equation. To this end, we introduce a Markovian framework, taking non-path-dependent coefficients: $b_t(x,a)$, $\sigma_t(x,a)$, $f_t(x,a)$ will depend on $(x,a)$ only through its value at time $t$, namely $(x(t),a(t))$. More precisely, in the present section we suppose that $b$, $\sigma$, $f$ are defined on $\R^n\times A$, with values respectively in $\R^n$, $\R^{n\times d}$, $\R$.

\subsection{Infinite horizon optimal control problem}

We consider the same probabilistic setting as in Section \ref{S:Formulation}, characterized by the following objects: $(\Omega,\Fc,\P)$, $W=(W_t)_{t\geq0}$, $\F^W=(\Fc_t^W)_{t\geq0}$, $\Ac$.

For every $x\in\R^n$ and $\alpha\in\Ac$, we consider the following Markovian controlled stochastic differential equation:
\begin{equation}\label{sde_Markov}
X_t^{x,\alpha} \ = \ x + \int_0^t b(X_s^{x,\alpha},\alpha_s)\,ds + \int_0^t \sigma(X_s^{x,\alpha},\alpha_s)\,dW_s,
\end{equation}
for all $t\geq0$. We then define the value function $v\colon\R^n\rightarrow\R$ of the corresponding infinite horizon stochastic optimal control problem as follows:
\begin{equation}\label{v}
v(x) \ = \ \sup_{\alpha\in\Ac} J(x,\alpha), \qquad \text{for all }x\in\R^n,
\end{equation}
where the gain functional is given by
\[
J(x,\alpha) \ = \ \E\bigg[\int_0^\infty e^{-\beta t} f(X_t^{x,\alpha},\alpha_t)\,dt\bigg].
\]
On the positive constant $\beta$ and on the coefficients $b$, $\sigma$, $f$ we impose the same sets of assumptions (either {\bf (A)} or {\bf (A')}) as for the path-dependent case. In the present Markovian framework those assumptions read as follows.

\vspace{3mm}

\noindent {\bf ($\mathbf A_{\text{\tiny Markov}}$)}
\begin{itemize}
\item [(i)]
The functions $b$, $\sigma$, $f$ are continuous.
\item [(iii)]
There exists a constant  $L\geq0$ such that
\begin{align*}
|b(x,a) - b(x',a)| + |\sigma(x,a)-\sigma(x',a)|
\ &\leq \ L\,|x-x'|,
\\
|b(0,a)| + |\sigma(0,a)| \ &\leq \ L,
\end{align*}
for all $x,x'\in\R^n$, $a\in A$.
\item[(iv)] There exist constants $M>0$ and $r\geq0$ such that
\[
|f(x,a)| \ \leq \ M(1 + |x|^r),
\]
for all $x\in\R^n$, $a\in A$.
\item[(v)] $\beta>\bar\beta$, where $\bar\beta$ is zero when the constant $r$ at point (iv) is zero; otherwise, when $r>0$, $\bar\beta$ is a strictly positive real number such that
\[
\E\Big[\sup_{s\in[0,t]}|X_s^{x,\alpha}|^r\Big] \ \leq \ \bar C e^{\bar\beta t}(1 + |x|^r),
\]
for some constant $\bar C\geq0$, with $\bar\beta$ and $\bar C$ independent of $t\geq0$, $\alpha\in\Ac$, $x\in\R^n$. See Lemma \ref{L:Estimate}.
\end{itemize}
\begin{Remark}
{\rm
Notice that, in the present Markovian framework, as the coefficients are independent of the time variable, we can unify {\bf (A)} and {\bf (A')} into a single set of assumptions {\bf ($\mathbf A_{\text{\tiny Markov}}$)}. As a matter of fact, in this case, if we include in {\bf (A')} the case $r=0$, with corresponding $\bar\beta=0$, then {\bf (A')} coincides with {\bf (A)}; while, in the path-dependent case, {\bf (A')} is in general stronger than {\bf (A)}, as the constant $L$ in {\bf (A')}-(iii) is taken independent of $T$.
\ep}
\end{Remark}
It is well-known that under {\bf ($\mathbf A_{\text{\tiny Markov}}$)}-(i)-(iii), for every $x\in\R^n$ and $\alpha\in\Ac$, there exists a unique $\F^W$-progressive continuous process $X^{x,\alpha}=(X_t^{x,\alpha})_{t\geq0}$ solution to equation \eqref{sde_Markov}. Moreover, under {\bf ($\mathbf A_{\text{\tiny Markov}}$)}-(i)-(iii), by Lemma \ref{L:Estimate}, we have, for every $p>0$,
\[
\E\Big[\sup_{s\in[0,t]}|X_s^{x,\alpha}|^p\Big] \ \leq \ \bar C_{p,L} e^{\bar\beta_{p,L} t}(1 + |x|^p),
\]
for some constants $\bar C_{p,L}\geq0$ and $\bar\beta_{p,L}\geq0$, independent of $t\geq0$, $\alpha\in\Ac$, $x\in\R^n$, depending only on $p>0$ and the constant $L$ appearing in {\bf ($\mathbf A_{\text{\tiny Markov}}$)}-(iii). Finally, we notice that the value function $v$ in \eqref{v} is well-defined by {\bf ($\mathbf A_{\text{\tiny Markov}}$)}-(iv)-(v).

\subsection{Randomized problem and backward SDE representation of $v(x)$}

We consider the same probabilistic setting as in Section \ref{S:Randomized}, with $\lambda$, $a_0$, $(\bar\Omega,\bar\Fc,\bar\P)$, $\bar W=(\bar W_t)_{t\geq0}$, $\bar\mu$ on $[0,\infty)\times A$ with compensator $\lambda(da)dt$, $\bar\F^{W,\mu}=(\bar\Fc_t^{W,\mu})_{t\geq0}$, $\bar I$ defined by \eqref{I}, $\bar\Vc$, and the family of probability measures $\bar\P_T^{\bar\nu}$, with $\bar\nu\in\bar\Vc$ and $T>0$. Then, for every $x\in\R^n$, we consider the stochastic differential equation
\begin{equation}\label{sde_random_Markov}
\bar X_t^x \ = \ x + \int_0^t b(\bar X_s^x,\bar I_s)\,ds + \int_0^t \sigma(\bar X_s^x,\bar I_s)\,d\bar W_s, \qquad \text{for all }t\geq0.
\end{equation}
It is well-known that under {\bf ($\mathbf A_{\text{\tiny Markov}}$)}-(i)-(iii), for every $x\in\R^n$, there exists a unique $\bar\F^{W,\mu}$-progressive continuous process $\bar X^x=(\bar X_t^x)_{t\geq0}$ solution to equation \eqref{sde_random_Markov}. Furthermore, under {\bf ($\mathbf A_{\text{\tiny Markov}}$)}-(i)-(iii), by Lemma \ref{L:Estimate}, for every $p>0$, there exist two non-negative constants $\bar C_{p,L}$ and $\bar\beta_{p,L}$ such that
\[
\bar\E_T^{\bar\nu}\Big[\sup_{s\in[0,T]}|\bar X_s^x|^p\Big|\bar\Fc_t^{W,\mu}\Big] \ \leq \ \bar C_{p,L} e^{\bar\beta_{p,L} t}\Big(1 + \sup_{s\in[0,t]}|\bar X_s^x|^p\Big),
\]
$\bar\P$-a.s., for all $t\in[0,T]$, $T>0$, $\bar\nu\in\bar\Vc$, where $\bar C_{p,L}$ and $\bar\beta_{p,L}$ are independent of $t\in[0,T]$, $T>0$, $\bar\nu\in\bar\Vc$, $x\in\R^n$, and depend only on $p>0$ and the constant $L$ appearing in {\bf ($\mathbf A_{\text{\tiny Markov}}$)}-(iii).

\begin{Remark}\label{R:Cbar_betabar}
{\rm
As we did in Lemma \ref{L:Estimate}, when $r>0$ in Assumption {\bf ($\mathbf A_{\text{\tiny Markov}}$)}\textup{-(iv)}, we denote $\bar C_{r,L}$ and $\bar\beta_{r,L}$ simply by $\bar C$ and $\bar\beta$.
\ep}
\end{Remark}

We define the value function $v^\Rc\colon\R^n\rightarrow\R$ of the randomized infinite horizon stochastic optimal control problem as follows:
\begin{equation}\label{value_R_Markov}
v^\Rc(x) \ = \ \lim_{T\rightarrow\infty}\sup_{\bar\nu\in\bar\Vc} J_T^\Rc(x,\bar\nu), \qquad \text{for all }x\in\R^n,
\end{equation}
where, for every $T>0$,
\[
J_T^\Rc(x,\bar\nu) \ = \ \bar\E_T^{\bar\nu}\bigg[\int_0^T e^{-\beta t} f(\bar X_t^x,\bar I_t)\,dt\bigg].
\]
From Section \ref{S:Randomized}, we know that the limit in \eqref{value_R_Markov} exists and it is finite for every $x\in\R^n$.

For every $x\in\R^n$, we now consider the following infinite horizon backward stochastic differential equation with nonpositive jumps:
\begin{align}
Y_t \ &= \ Y_T - \beta\int_t^T Y_s\,ds + \int_t^T f(\bar X_s^x,\bar I_s)\,ds + K_T - K_t \notag \\
&\quad \ - \int_t^T Z_s\,d\bar W_s - \int_t^T \int_A U_s(a)\,\bar\mu(ds,da), \qquad 0 \leq t \leq T, \;  \forall\,T \in (0,\infty), \; \bar\P\text{-a.s.} \label{BSDE_Markov} \\
\notag \\
U_t(a) \ &\leq \ 0, \qquad dt\otimes d\bar\P\otimes\lambda(da)\text{-a.e. on }[0,\infty)\times\bar\Omega\times A. \label{JumpConstr_Markov}
\end{align}

\begin{Lemma}\label{L:BSDE_x}
Under {\bf ($\mathbf A_{\text{\tiny Markov}}$)}, for every $x\in\R^n$, there exists a unique quadruple $(\bar Y^x,\bar Z^x,\bar U^x,\bar K^x)$ in ${\bf S_{\textup{loc}}^2}\times{\bf L_{\textup{loc}}^2(W)}\times{\bf L_{\textup{loc}}^2(\bar\mu)}\times{\bf K_{\textup{loc}}^2}$, such that $\bar Y^x\in{\bf S^\infty}$ $($resp. $|\bar Y_t^x|\leq C(1+\sup_{s\in[0,t]}|\bar X_s^x|^r)$, for all $t\geq0$, $\bar\P$-a.s., and for some positive constant $C$$)$, satisfying \eqref{BSDE_Markov}-\eqref{JumpConstr_Markov} which is minimal in the following sense: for any other quadruple $(\hat Y^x,\hat Z^x,\hat U^x,\hat K^x)$ in ${\bf S_{\textup{loc}}^2}\times{\bf L_{\textup{loc}}^2(W)}\times{\bf L_{\textup{loc}}^2(\bar\mu)}\times{\bf K_{\textup{loc}}^2}$, with $\hat Y^x\in{\bf S^\infty}$ $($resp. $|\hat Y_t^x|\leq C(1+\sup_{s\in[0,t]}|\bar X_s^x|^r)$, for all $t\geq0$, $\bar\P$-a.s., and for some positive constant $C$$)$, satisfying \eqref{BSDE_Markov}-\eqref{JumpConstr_Markov} we have:
\[
\bar Y_t^x \ \leq \ \hat Y_t^x, \qquad \text{for all }t\geq0, \; \bar\P\text{-a.s.}
\]
\end{Lemma}
\textbf{Proof.}
The result follows directly from Theorem \ref{T:Main}.
\ep

\begin{Proposition}\label{P:BSDE_Markov}
Suppose that {\bf ($\mathbf A_{\text{\tiny Markov}}$)} holds. Then, for every $x\in\R^n$, we have:
\begin{itemize}
\item[\textup{(i)}] $v(x)=v^\Rc(x)$;
\item[\textup{(ii)}] $v(x)=\bar Y_0^x$, $\bar\P$-a.s., where $(\bar Y^x,\bar Z^x,\bar U^x,\bar K^x)$ is the minimal solution to \eqref{BSDE_Markov}-\eqref{JumpConstr_Markov}; moreover, we have $($$r$ is the constant appearing in {\bf ($\mathbf A_{\text{\tiny Markov}}$)}\textup{-(iv)}$)$:
\begin{itemize}
\item[$r=0:$] $|v(x)|\leq \frac{\|f\|_\infty}{\beta}$, for all $x\in\R^n$;
\item[$r>0:$] $|v(x)|\leq \frac{2M(1+\bar C)}{\beta-\bar\beta}(1 + |x|^r)$, for all $x\in\R^n$;
\end{itemize}
\item[\textup{(iii)}] $v(\bar X_t^x)=\bar Y_t^x$, $\bar\P$-a.s., for all $t\geq0$;
\item[\textup{(iv)}] $v$ satisfies the following equality:
\begin{equation}\label{DPP_Markov}
v(x) \ = \ \sup_{\bar\nu\in\bar\Vc}\bar\E_T^{\bar\nu}\bigg[\int_0^\tau e^{-\beta s}\,f(\bar X_s^x,\bar I_s)\,ds + e^{-\beta \tau}\,v(\bar X_\tau^x)\bigg],
\end{equation}
for all $T\geq0$ and any $\bar\F^{W,\mu}$-stopping time $\tau$ taking values in $[0,T]$.
\end{itemize}
\end{Proposition}

\begin{Remark}
{\rm
We refer to formula \eqref{DPP_Markov} as the \emph{randomized dynamic programming principle} for $v$.
\ep}
\end{Remark}
\textbf{Proof (of Proposition \ref{P:BSDE_Markov}).}
Point (i) is a consequence of Theorem \ref{T:Identification}, while point (ii) follows from Theorem \ref{T:Main}, Lemma \ref{L:BSDE_x}, and either estimate \eqref{bound_Y} (case $r=0$) or estimate \eqref{bound_Y_pol} (case $r>0$) with $t=0$. Concerning identity (iii), for every $n\in\N$ consider the penalized infinite horizon backward stochastic differential equation \eqref{BSDE_n} with $f_s(\bar X,\bar I)$ replaced by $f(\bar X_s^x,\bar I_s)$, and denote by $(\bar Y^{n,x},\bar Z^{n,x},\bar U^{n,x})\in{\bf S_{\textup{loc}}^2}\times{\bf L_{\textup{loc}}^2(W)}\times{\bf L_{\textup{loc}}^2(\bar\mu)}$ its solution. Then, it is a standard result in the theory of backward stochastic differential equations that, for every $n$, there exists a function $v_n\colon\R^n\rightarrow\R$, with $v_n(x)=\bar Y_0^{n,x}$, satisfying $v_n(\bar X_t^x)=\bar Y_t^{n,x}$, $\bar\P$-a.s., for all $t\geq0$. By step II of the proof of Theorem \ref{T:Main}, we know that $\bar Y_0^{n,x}=v_n(x)$ converges increasingly to $\bar Y_0^x=v(x)$. Then, letting $n\rightarrow\infty$ in identity $v_n(\bar X_t^x)=\bar Y_t^{n,x}$, we deduce (iii). Finally, formula \eqref{DPP_Markov} follows from point (iii) and formula \eqref{DPP}, with $t=0$.
\ep

\subsection{Fully non-linear elliptic PDE}

Consider the following Hamilton-Jacobi-Bellman fully non-linear elliptic partial differential equation on $\R^n$:
\begin{equation}\label{HJB}
\beta\,u(x) - \sup_{a\in A}\big[\Lc^a u(x) + f(x,a)\big] \ = \ 0, \qquad \text{for all }x\in\R^n,
\end{equation}
where
\[
\Lc^a u(x) \ = \ \langle b(x,a),D_x u(x)\rangle + \frac{1}{2}\text{tr}\big[\sigma\sigma\trans(x,a)D_x^2 u(x)\big].
\]
We now prove, by means of the randomized dynamic programming principle,
that the value function $v$ in \eqref{v} is a viscosity solution of equation \eqref{HJB}.
As explained in the introduction, this allows us to avoid
 the difficulties related to a rigorous
proof of the classical dynamic programming principle,
which often  requires the use of nontrivial
measurability arguments.
In order to do this, we need an additional assumption.

\vspace{3mm}

\noindent {\bf ($\mathbf A_f$)} The function $f=f(x,a)\colon\R^n\times A\rightarrow\R$ is continuous in $x$ uniformly with respect to $a$.

\begin{Remark}\label{R:ContH}
{\rm
Notice that under {\bf ($\mathbf A_{\text{\tiny Markov}}$)} and {\bf ($\mathbf A_f$)}, for every $\varphi\in C^2(\R^n)$, the map
\[
x\in\R^n \ \longmapsto \ \sup_{a\in A}\big[\Lc^a \varphi(x) + f(x,a)\big]
\]
is continuous on $\R^n$.
\ep}
\end{Remark}

Given a locally bounded function $u\colon\R^n\rightarrow\R$, we define its upper and lower semi-continuous envelopes:
\[
u^*(x) \ = \ \limsup_{y\rightarrow x}u(y), \qquad\qquad u_*(x) \ = \ \liminf_{y\rightarrow x}u(y),
\]
for all $x\in\R^n$.

\begin{Theorem}\label{T:Markov}
Under {\bf ($\mathbf A_{\text{\tiny Markov}}$)} and {\bf ($\mathbf A_f$)}, the value function $v$ in \eqref{v} is a viscosity solution to \eqref{HJB}, namely $v$ is locally bounded and satisfies:
\begin{itemize}
\item $v^*$ is a viscosity subsolution to \eqref{HJB}:
\[
\beta\,\varphi(x) - \sup_{a\in A}\big[\Lc^a \varphi(x) + f(x,a)\big] \ \leq \ 0,
\]
for every $x\in\R^n$ and $\varphi\in C^2(\R^n)$ such that
\begin{equation}\label{ViscSub}
\max_{y\in\R^n} \big(v^*(y) - \varphi(y)\big) \ = \ v^*(x) - \varphi(x) \ = \ 0.
\end{equation}
\item $v_*$ is a viscosity supersolution to \eqref{HJB}:
\[
\beta\,\varphi(x) - \sup_{a\in A}\big[\Lc^a \varphi(x) + f(x,a)\big] \ \geq \ 0,
\]
for every $x\in\R^n$ and $\varphi\in C^2(\R^n)$ such that
\begin{equation}\label{ViscSuper}
\min_{y\in\R^n} \big(v_*(y) - \varphi(y)\big) \ = \ v_*(x) - \varphi(x) \ = \ 0.
\end{equation}
\end{itemize}
\end{Theorem}
\begin{Remark}\label{HJBunique}
{\rm
We do not provide here a comparison principle (and hence a uniqueness result) for the fully non-linear elliptic partial differential equation \eqref{HJB}, which can be found for instance in \cite{Ishii}, Theorem 7.4, under stronger assumptions than {\bf ($\mathbf A_{\text{\tiny Markov}}$)}-{\bf ($\mathbf A_f$)}. However, we notice that, when a uniqueness result for equation \eqref{HJB} holds, Theorem \ref{T:Markov} provides a non-linear Feynman-Kac representation formula for the unique viscosity solution (reported in Proposition \ref{P:BSDE_Markov}-(ii)), expressed in terms of the class of infinite horizon backward stochastic differential equations with nonpositive jumps \eqref{BSDE_Markov}-\eqref{JumpConstr_Markov}, with $x\in\R^n$.
\ep}
\end{Remark}
\textbf{Proof.} We begin noting that $v$ is locally bounded, as a consequence of point (ii) in Proposition \ref{P:BSDE_Markov}. Now, we split the proof into two steps.

\vspace{2mm}

\noindent\textbf{Step I.} \emph{Viscosity subsolution property.} Take $x\in\R^n$ and $\varphi\in C^2(\R^n)$ such that \eqref{ViscSub} holds. Without loss of generality, we suppose that the maximum in \eqref{ViscSub} is strict. As a matter of fact, if we prove the viscosity subsolution property for all $\varphi$ for which the maximum in \eqref{ViscSub} is strict, then given a generic $\varphi\in C^2(\R^n)$ satisfying \eqref{ViscSub}, it is enough to consider the function $\psi\in C^2(\R^n)$ given by $\psi(y)=\varphi(y)+|y-x|^4$, for all $y\in\R^n$, and to notice that $\psi(x)=\varphi(x)$, $D_x\psi(x)=D_x\varphi(x)$, $D_x^2\psi(x)=D_x^2\varphi(x)$; moreover, $\psi$ satisfies \eqref{ViscSub} and the maximum is strict.

We proceed by contradiction, assuming that
\[
\beta\,\varphi(x) - \sup_{a\in A}\big[\Lc^a \varphi(x) + f(x,a)\big] \ =: \ 2\,\eps_0 \ > \ 0.
\]
By Remark \ref{R:ContH}, it follows that the map $y\mapsto\beta\,\varphi(y) - \sup_{a\in A}[\Lc^a \varphi(y) + f(y,a)]$ is continuous on $\R^n$, therefore there exists $\eta>0$ such that
\[
\beta\,\varphi(y) - \sup_{a\in A}\big[\Lc^a \varphi(y) + f(y,a)\big] \ \geq \ \eps_0, \qquad \text{for all }y\in\R^n\text{ satisfying }|y-x| \ \leq \ \eta.
\]
Since the maximum in \eqref{ViscSub} is strict, we have
\begin{equation}\label{delta}
\max_{y\in\R^n\,\colon|y-x|=\eta} \big(v^*(y) - \varphi(y)\big) \ = \ -\delta \ < \ 0.
\end{equation}
Fix $\eps\in(0,\min(\eps_0,\delta\beta)]$. From the definition of $v^*(x)$, there exists a sequence $\{y_m\}_{m\in\N}\subset\R^n$, with $|y_m-x|<\eta$, such that
\[
y_m \ \overset{m\rightarrow\infty}{\longrightarrow} \ x, \qquad\qquad v(y_m) \ \overset{m\rightarrow\infty}{\longrightarrow} \ v^*(x).
\]
Since $\varphi(y_m)\rightarrow\varphi(x)$ as $m\rightarrow\infty$, and $\varphi(x)=v^*(x)$, we have
\[
\gamma_m \ := \ v(y_m) - \varphi(y_m) \ \overset{m\rightarrow\infty}{\longrightarrow} \ 0.
\]
Take $T>0$ and set
\[
\tau_m \ := \ \inf\big\{t\geq0\colon\big|\bar X_t^{y_m} - x\big| \ \geq \ \eta\big\}, \qquad\qquad \theta_m \ := \ \tau_m \ \wedge \ T,
\]
where $\bar X^{y_m}=(\bar X_t^{y_m})_{t\geq0}$ denotes the solution to equation \eqref{sde_random_Markov} with $x$ replaced by $y_m$. Notice that $\tau_m$ and $\theta_m$ are $\bar\F^{W,\mu}$-stopping times, with $\theta_m$ taking values in $[0,T]$. Then, by formula \eqref{DPP_Markov}, it follows that there exists $\bar\nu^m\in\bar\Vc$ such that
\[
v(y_m) - \frac{\eps}{2\beta} \ \leq \ \bar\E_T^{\bar\nu^m}\bigg[\int_0^{\theta_m} e^{-\beta s}\,f(\bar X_s^{y_m},\bar I_s)\,ds + e^{-\beta \theta_m}\,v(\bar X_{\theta_m}^{y_m})\bigg].
\]
By $v\leq v^*\leq\varphi$ and \eqref{delta}, we deduce that
\[
\gamma_m - \frac{\eps}{2\beta} \ \leq \ \bar\E_T^{\bar\nu^m}\bigg[\int_0^{\theta_m} e^{-\beta s}\,f(\bar X_s^{y_m},\bar I_s)\,ds + e^{-\beta \theta_m}\,\varphi(\bar X_{\theta_m}^{y_m}) - \varphi(y_m) - \delta\,e^{-\beta \theta_m}\,1_{\{\tau_m\leq T\}}\bigg].
\]
An application of It\^o's formula to $e^{-\beta s}\,\varphi(\bar X_s^{y_m})$ between $s=0$ and $s=\theta_m$, yields
\[
\gamma_m - \frac{\eps}{2\beta} \ \leq \ -\,\bar\E_T^{\bar\nu^m}\bigg[\int_0^{\theta_m} e^{-\beta s}\,\big(\beta\,\varphi(\bar X_s^{y_m}) - \Lc^{\bar I_s}\varphi(\bar X_s^{y_m}) - f(\bar X_s^{y_m},\bar I_s)\big)\,ds + \delta\,e^{-\beta \theta_m}\,1_{\{\tau_m\leq T\}}\bigg].
\]
Now, we observe that, whenever $0\leq s\leq\theta_m$,
\[
\beta\,\varphi(\bar X_s^{y_m}) - \Lc^{\bar I_s}\varphi(\bar X_s^{y_m}) - f(\bar X_s^{y_m},\bar I_s) \ \geq \ \beta\,\varphi(\bar X_s^{y_m}) - \sup_{a\in A}\big[\Lc^a\varphi(\bar X_s^{y_m}) + f(\bar X_s^{y_m},a)\big] \ \geq \ \eps.
\]
Therefore
\begin{align*}
\gamma_m - \frac{\eps}{2\beta} \ &\leq \ -\,\bar\E_T^{\bar\nu^m}\bigg[\int_0^{\theta_m} \eps\,e^{-\beta s}\,ds + \delta\,e^{-\beta \theta_m}\,1_{\{\tau_m\leq T\}}\bigg] \\
&= \ - \frac{\eps}{\beta} + \bar\E_T^{\bar\nu^m}\bigg[\Big(\frac{\eps}{\beta} - \delta\Big)\,e^{-\beta \theta_m}\,1_{\{\tau_m\leq T\}} + \frac{\eps}{\beta}\,e^{-\beta \theta_m}\,1_{\{\tau_m>T\}}\bigg] \\
&\leq \ - \frac{\eps}{\beta} + \frac{\eps}{\beta}\,\bar\E_T^{\bar\nu^m}\big[e^{-\beta \theta_m}\,1_{\{\tau_m>T\}}\big] \ = \ - \frac{\eps}{\beta} + \frac{\eps}{\beta}\,\bar\E_T^{\bar\nu^m}\big[e^{-\beta T}\,1_{\{\tau_m>T\}}\big] \\
&= \ - \frac{\eps}{\beta} + \frac{\eps}{\beta}\,e^{-\beta T}\,\bar\P_T^{\bar\nu^m}(\tau_m \ > \ T) \ \leq \ - \frac{\eps}{\beta} + \frac{\eps}{\beta}\,e^{-\beta T}.
\end{align*}
Letting $T\rightarrow\infty$, we conclude that
\[
\gamma_m - \frac{\eps}{2\beta} \ \leq \ - \frac{\eps}{\beta},
\]
Sending $m\rightarrow\infty$, we find the contradiction $-\eps/(2\beta)\leq-\eps/\beta$.


\vspace{2mm}

\noindent\textbf{Step II.} \emph{Viscosity supersolution property.} Take $x\in\R^n$ and $\varphi\in C^2(\R^n)$ such that \eqref{ViscSuper} holds. From the definition of $v_*(x)$, there exists a sequence $\{y_m\}_{m\in\N}\subset\R^n$ such that
\[
y_m \ \overset{m\rightarrow\infty}{\longrightarrow} \ x, \qquad\qquad v(y_m) \ \overset{m\rightarrow\infty}{\longrightarrow} \ v_*(x).
\]
Since $\varphi(y_m)\rightarrow\varphi(x)$ as $m\rightarrow\infty$, and $\varphi(x)=v_*(x)$, we also have
\[
\gamma_m \ := \ v(y_m) - \varphi(y_m) \ \overset{m\rightarrow\infty}{\longrightarrow} \ 0.
\]
Consider a sequence $\{h_m\}_{m\in\N}\subset(0,\infty)$ such that
\[
h_m \ \overset{m\rightarrow\infty}{\longrightarrow} \ 0, \qquad\qquad \frac{\gamma_m}{h_m} \ \overset{m\rightarrow\infty}{\longrightarrow} \ 0.
\]
Let us denote by $\bar X^{y_m}=(\bar X_t^{y_m})_{t\geq0}$ the solution to equation \eqref{sde_random_Markov} with $x$ replaced by $y_m$. Given $\eta>0$, we define
\[
\tau_m \ := \ \inf\big\{t\geq0\colon\big|\bar X_t^{y_m} - x\big| \ \geq \ \eta\big\}, \qquad\qquad \theta_m \ := \ \tau_m \ \wedge \ h_m.
\]
Then, by formula \eqref{DPP_Markov}, we obtain (recall that the map $\bar\nu\equiv1$ belongs to $\bar\Vc$ and, for such $\bar\nu$, we have that $\bar\P_T^{\bar\nu}$ is $\bar\P$, for any $T>0$)
\[
v(y_m) \ \geq \ \bar\E\bigg[\int_0^{\theta_m} e^{-\beta s}\,f(\bar X_s^{y_m},\bar I_s)\,ds + e^{-\beta \theta_m}\,v(\bar X_{\theta_m}^{y_m})\bigg].
\]
Since $v\geq v_*\geq\varphi$, we obtain
\[
\gamma_m \ \geq \ \bar\E\bigg[\int_0^{\theta_m} e^{-\beta s}\,f(\bar X_s^{y_m},\bar I_s)\,ds + e^{-\beta \theta_m}\,\varphi(\bar X_{\theta_m}^{y_m}) - \varphi(y_m)\bigg].
\]
Applying It\^o's formula to $e^{-\beta s}\varphi(\bar X_s^{y_m})$ between $s=0$ and $s=\theta_m$, we find
\begin{equation}\label{gamma_h}
\frac{\gamma_m}{h_m} \ \geq \ -\,\bar\E\bigg[\frac{1}{h_m}\int_0^{\theta_m} e^{-\beta s}\,\big(\beta\,\varphi(\bar X_s^{y_m}) - \Lc^{\bar I_s}\varphi(\bar X_s^{y_m}) - f(\bar X_s^{y_m},\bar I_s)\big)\,ds\bigg].
\end{equation}
Now, notice that, for $\bar\P$-a.e. $\bar\omega\in\bar\Omega$, there exists
$M(\bar\omega)\in\N$ such that $\theta_m(\bar\omega)=h_m$ for any
$m\geq M(\bar\omega)$. Moreover, $\bar I_s$ converges $\bar\P$-a.s. to $a_0$ as $s\searrow 0$.
In addition, we recall the following standard result:
\[
\bar\E\Big[\sup_{|y-x|\leq\eta}\sup_{s\in[0,S]}\big|\bar X_s^y - x\big|\Big] \ \overset{S\rightarrow0}{\longrightarrow} \ 0.
\]
Therefore, up to a subsequence, we have
\[
\xi_m := \frac{1}{h_m}\int_0^{\theta_m} e^{-\beta s}\,\big(\beta\,\varphi(\bar X_s^{y_m}) - \Lc^{\bar I_s}\varphi(\bar X_s^{y_m}) - f(\bar X_s^{y_m},\bar I_s)\big)\,ds \ \underset{\bar\P\text{-a.s.}}{\overset{m\rightarrow\infty}{\longrightarrow}} \ \beta\,\varphi(x) - \Lc^{a_0} \varphi(x) - f(x,a_0).
\]
Furthermore, since the integral in $\xi_m$ is over the interval $[0,\theta_m]$, it follows that there exists a non-negative $\xi\in L^1(\bar\Omega,\bar\Fc,\bar\P)$ such that $|\xi_m|\leq\xi$, $\bar\P$-a.s., for every $m\in\N$. Then, by the Lebesgue dominated convergence theorem, sending $m\rightarrow\infty$ in \eqref{gamma_h}, we find
\[
\beta\,\varphi(x) \ \geq \ \Lc^{a_0} \varphi(x) + f(x,a_0).
\]
Recalling that the deterministic point $a_0\in A$ fixed at the beginning of Section \ref{S:Randomized} was arbitrary, we conclude that
\[
\beta\,\varphi(x) \ \geq \ \sup_{a\in A}\big[\Lc^a \varphi(x) + f(x,a)\big].
\]
\ep

\appendix

\setcounter{equation}{0} \setcounter{Assumption}{0}
\setcounter{Theorem}{0} \setcounter{Proposition}{0}
\setcounter{Corollary}{0} \setcounter{Lemma}{0}
\setcounter{Definition}{0} \setcounter{Remark}{0}

\renewcommand\thesection{Appendix}

\section{}

\renewcommand\thesection{\Alph{subsection}}

\renewcommand\thesubsection{\hspace{-5mm}}

\subsection{Proof of Lemma \ref{L:Estimate}}

\textbf{Proof.}
We begin proving estimate \eqref{beta_hat} for $p\geq2$. Given $t\geq0$, $A\in\hat\Fc_t$, $T\geq t$, we obtain, from equation \eqref{sde_hat},
\begin{align}
\hat\E\Big[\sup_{s\in[0,T]}|\hat X_s|^p\,1_A\Big]^{1/p} \ &\leq \ \hat\E\Big[\sup_{s\in[0,t]}|\hat X_s|^p\,1_A\Big]^{1/p} + \hat\E\bigg[\sup_{s\in[t,T]}\bigg|\int_t^s b_u(\hat X,\gamma)\,du\bigg|^p\,1_A\bigg]^{1/p} \notag \\
&\quad \ +
\hat\E\bigg[\sup_{s\in[t,T]}\bigg|\int_t^s \sigma_u(\hat X,\gamma)\,d\hat W_u\bigg|^p\,1_A\bigg]^{1/p}. \label{estimate_proof1}
\end{align}
Notice that
\[
\hat\E\bigg[\sup_{s\in[t,T]}\bigg|\int_t^s b_u(\hat X,\gamma)\,du\bigg|^p\,1_A\bigg]^{1/p} \ \leq \ \int_t^T \hat\E\big[\big|b_u(\hat X,\gamma)\big|^p\,1_A\big]^{1/p}\,du.
\]
Then, by Assumption {\bf (A)'}-(iii), we obtain
\begin{align}
\hat\E\bigg[\sup_{s\in[t,T]}\bigg|\int_t^s b_u(\hat X,\gamma)\,du\bigg|^p\,1_A\bigg]^{1/p} \ &\leq \ L\int_t^T \hat\E\Big[\Big(1+\sup_{s\in[0,u]}|\hat X_s|\Big)^p\,1_A\Big]^{1/p}\,du \notag \\
&\leq \ L\,(T-t)\,\hat\E[1_A]^{1/p} + L\int_t^T \hat\E\Big[\sup_{s\in[0,u]}|\hat X_s|^p\,1_A\Big]^{1/p}\,du. \label{b_estimate}
\end{align}
On the other hand, by the Burkholder-Davis-Gundy inequality,
\begin{align*}
\hat\E\bigg[\sup_{s\in[t,T]}\bigg|\int_t^s \sigma_u(\hat X,\gamma)\,d\hat W_u\bigg|^p\,1_A\bigg]^{1/p} \ &\leq \ C_p\,\bigg(\hat\E\bigg[\bigg(\int_t^T |\sigma_u(\hat X,\gamma)|^2\,du\bigg)^{p/2}\,1_A\bigg]\bigg)^{1/p} \\
&= \ C_p\,\bigg\|\int_t^T |\sigma_u(\hat X,\gamma)|^2\,1_A\,du\bigg\|_{L^{p/2}(\hat\Omega,\hat\Fc,\hat\P)}^{1/2} \\
&\leq \ C_p\,\bigg(\int_t^T \hat\E\big[|\sigma_u(\hat X,\gamma)|^p\,1_A\big]^{2/p}\,du\bigg)^{1/2},
\end{align*}
where $C_p>0$ is the constant of the Burkholder-Davis-Gundy inequality. By Assumption {\bf (A)'}-(iii), we find
\begin{align*}
\hat\E\bigg[\sup_{s\in[t,T]}\bigg|\int_t^s \sigma_u(\hat X,\gamma)\,d\hat W_u\bigg|^p\,1_A\bigg]^{1/p} \ \leq \ C_p\,L\,\bigg(\int_t^T \hat\E\Big[\Big(1+\sup_{s\in[0,u]}|\hat X_s|\Big)^p\,1_A\Big]^{2/p}\,du\bigg)^{1/2}& \\
\leq \ C_p\,L\,(T-t)\,\hat\E[1_A]^{1/p} + C_p\,L\bigg(\int_t^T \hat\E\Big[\sup_{s\in[0,u]}|\hat X_s|^p\,1_A\Big]^{2/p}\,du\bigg)^{1/2}&.
\end{align*}
Plugging this latter estimate and \eqref{b_estimate} into \eqref{estimate_proof1}, we obtain
\begin{align}\label{estimate_proof2}
\hat\E\Big[\sup_{s\in[0,T]}|\hat X_s|^p\,1_A\Big]^{1/p} \ &\leq \ \hat\E\Big[\sup_{s\in[0,t]}|\hat X_s|^p\,1_A\Big]^{1/p} + L\int_t^T \hat\E\Big[\sup_{s\in[0,u]}|\hat X_s|^p\,1_A\Big]^{1/p}\,du \\
&\quad \ +
(1 + C_p)\,L\,(T-t)\,\hat\E[1_A]^{1/p} + C_p\,L\bigg(\int_t^T \hat\E\Big[\sup_{s\in[0,u]}|\hat X_s|^p\,1_A\Big]^{2/p}\,du\bigg)^{1/2}. \notag
\end{align}
Taking the square of both sides of \eqref{estimate_proof2}, and using the inequality $(x_1+x_2+x_3+x_4)^2\leq 4(x_1^2+x_2^2+x_3^2+x_4^2)$, valid for any $x_1,x_2,x_3,x_4\in\R$, we find
\begin{align}\label{estimate_proof3}
\hat\E\Big[\sup_{s\in[0,T]}|\hat X_s|^p\,1_A\Big]^{2/p} &\leq 4\,\hat\E\Big[\sup_{s\in[0,t]}|\hat X_s|^p\,1_A\Big]^{2/p} + 4\,L^2\bigg(\int_t^T \hat\E\Big[\sup_{s\in[0,u]}|\hat X_s|^p\,1_A\Big]^{1/p}\,du\bigg)^2 \\
&\quad \ +
4\,(1 + C_p)^2\,L^2\,(T-t)^2\,\hat\E[1_A]^{2/p} + 4\,C_p^2\,L^2\int_t^T \!\!\hat\E\Big[\sup_{s\in[0,u]}|\hat X_s|^p\,1_A\Big]^{2/p}du. \notag
\end{align}
Now, consider the nonnegative function $v\colon[t,\infty)\rightarrow[0,\infty)$ given by
\[
v(T) \ := \ \hat\E\Big[\sup_{s\in[0,T]}|\hat X_s|^p\,1_A\Big]^{2/p}, \qquad \text{for all }T\geq t.
\]
Then, rewriting inequality \eqref{estimate_proof3} in terms of $v$, we obtain
\begin{align}\label{estimate_proof4}
v(T) \ &\leq \ 4\,\hat\E\Big[\sup_{s\in[0,t]}|\hat X_s|^p\,1_A\Big]^{2/p} + 4\,(1 + C_p)^2\,L^2\,(T-t)^2\,\hat\E[1_A]^{2/p} + 4\,L^2\bigg(\int_t^T \sqrt{v(u)}\,du\bigg)^2 \notag \\
&\quad \ + 4\,C_p^2\,L^2\,\int_t^T v(u)\,du.
\end{align}
Multiplying both sides of \eqref{estimate_proof4} by $e^{-(T-t)}$:
\begin{align}\label{estimate_proof5}
e^{-(T-t)}\,v(T) \ &\leq \ 4\,\hat\E\Big[\sup_{s\in[0,t]}|\hat X_s|^p\,1_A\Big]^{2/p} + 4\,(1 + C_p)^2\,L^2\,(T-t)^2\,\hat\E[1_A]^{2/p} \notag \\
&\quad \ + 4\,L^2\bigg(\int_t^T e^{-(T-t)/2}\,\sqrt{v(u)}\,du\bigg)^2 + 4\,C_p^2\,L^2\,\int_t^T e^{-(T-t)}v(u)\,du,
\end{align}
where we have used the inequality $e^{-(T-t)}\leq1$. Now, notice that
\begin{align*}
\bigg(\int_t^T e^{-(T-t)/2}\,\sqrt{v(u)}\,du\bigg)^2 \ &= \ \bigg(\int_t^T e^{-(T-u)/2}\,\big(e^{-(u-t)}\,v(u)\big)^{1/2}\,du\bigg)^2 \\
&\leq \ \int_t^T e^{-(T-u)}\,du \int_t^T e^{-(u-t)}\,v(u)\,du \ \leq \ \int_t^T e^{-(u-t)}\,v(u)\,du
\end{align*}
and
\[
\int_t^T e^{-(T-t)}v(u)\,du \ \leq \ \int_t^T e^{-(u-t)}v(u)\,du.
\]
Therefore, from inequality \eqref{estimate_proof5}, we get
\begin{align*}
e^{-(T-t)}\,v(T) \ &\leq \ 4\,\hat\E\Big[\sup_{s\in[0,t]}|\hat X_s|^p\,1_A\Big]^{2/p} + 4\,(1+ C_p)^2\,L^2\,(T-t)^2\,\hat\E[1_A]^{2/p} \\
&\quad \ + 4\,(1 + C_p^2)\,L^2\,\int_t^T e^{-(u-t)}\,v(u)\,du.
\end{align*}
An application of Gronwall's lemma to the function $u\mapsto e^{-(u-t)}v(u)$, yields
\[
e^{-(T-t)}\,v(T) \ \leq \ 4\,\Big(\hat\E\Big[\sup_{s\in[0,t]}|\hat X_s|^p\,1_A\Big]^{2/p} + (1 + C_p)^2\,L^2\,(T-t)^2\,\hat\E[1_A]^{2/p}\Big)\,e^{4\,(1 + C_p^2)\,L^2\,(T-t)}.
\]
Then, noting that $(T-t)^2\leq e^{T-t}$, we end up with
\[
\hat\E\Big[\sup_{s\in[0,T]}|\hat X_s|^p\,1_A\Big]^{2/p} \ \leq \ 4\,\Big(\hat\E\Big[\sup_{s\in[0,t]}|\hat X_s|^p\,1_A\Big]^{2/p} + (1 + C_p)^2\,L^2\,\hat\E[1_A]^{2/p}\Big)\,e^{(1 + 4\,(1 + C_p^2)\,L^2)\,(T-t)},
\]
from which we find (using the inequality $(a+b)^{p/2}\leq 2^{p/2-1}(a^{p/2}+b^{p/2})$, valid for any $a,b\geq0$)
\begin{align*}
\hat\E\Big[\sup_{s\in[0,T]}|\hat X_s|^p\,1_A\Big] \ &\leq \ 2^{p+p/2-1}\,\Big(\hat\E\Big[\sup_{s\in[0,t]}|\hat X_s|^p\,1_A\Big] + (1+ C_p)^p\,L^p\,\hat\E[1_A]\Big)\,e^{\frac{p}{2}(1 + 4\,(1 + C_p^2)\,L^2)\,(T-t)} \\
&\leq \ \bar C_{p,L}\,\hat\E\Big[\Big(1 + \sup_{s\in[0,t]}|\hat X_s|^p\Big)\,1_A\Big]\,e^{\bar\beta_{p,L}(T-t)},
\end{align*}
with (for $p\geq2$)
\[
\bar C_{p,L} \ := \ 2^{p+p/2-1}\max\big(1,(1+ C_p)^pL^p\big), \qquad\qquad \bar\beta_{p,L} \ := \ \frac{p}{2}\big(1 + 4(1 + C_p^2)L^2\big).
\]
From the arbitrariness of $A\in\hat\Fc_t$, we conclude
\begin{equation}\label{estimate_proof6}
\hat\E\Big[\sup_{s\in[0,T]}|\hat X_s|^p\Big|\hat\Fc_t\Big] \ \leq \ \bar C_{p,L}\,\Big(1 + \sup_{s\in[0,t]}|\hat X_s|^p\Big)\,e^{\bar\beta_{p,L}(T-t)}, \qquad \hat\P\text{-a.s.}
\end{equation}
which yields estimate \eqref{beta_hat} for $p\geq2$

Finally, when $p\in(0,2)$, we have, by Jensen's inequality,
\[
\hat\E\Big[\sup_{s\in[0,T]}|\hat X_s|^p\Big|\hat\Fc_t\Big]^{2/p} \ \leq \ \hat\E\Big[\sup_{s\in[0,T]}|\hat X_s|^2\Big|\hat\Fc_t\Big].
\]
Then, by estimate \eqref{estimate_proof6} with $p=2$, we obtain
\begin{align*}
\hat\E\Big[\sup_{s\in[0,T]}|\hat X_s|^p\Big|\hat\Fc_t\Big] \ &\leq \ \bar C_{2,L}^{p/2}\,\Big(1 + \sup_{s\in[0,t]}|\hat X_s|^2\Big)^{p/2}\,e^{\frac{p}{2}\bar\beta_{2,L}(T-t)} \\
&\leq \ \bar C_{2,L}^{p/2}\,\Big(1 + \sup_{s\in[0,t]}|\hat X_s|^p\Big)\,e^{\frac{p}{2}\bar\beta_{2,L}(T-t)},
\end{align*}
which yields estimate \eqref{beta_hat} for $p\in(0,2)$, with
\[
\bar C_{p,L} \ := \ \bar C_{2,L}^{p/2} \ = \ 2^p\max\big(1,(1+ C_2)^pL^p\big), \qquad \bar\beta_{p,L} \ := \ \frac{p}{2}\bar\beta_{2,L} \ = \ \frac{p}{2}\big(1 + 4(1 + C_2^2)L^2\big).
\]
\ep

\end{document}